\documentclass[12pt,a4paper, leqno]{article}
\usepackage{amsmath}
\usepackage[active]{srcltx}
\usepackage{amssymb}
\usepackage[T1]{fontenc}
\usepackage{url}
\usepackage{textcomp}
\usepackage{amsthm}
\usepackage{dsfont}
\usepackage{mathrsfs}
\usepackage[dvips]{graphicx}
\usepackage{psfrag}
\usepackage[hang,small,bf]{caption}
\usepackage{fullpage}
\newtheorem{theorem}{Theorem}[section]
\newtheorem{corollary}[theorem]{Corollary}
\newtheorem{definition}[theorem]{Definition}
\newtheorem{proposition}[theorem]{Proposition}
\newtheorem{lemma}[theorem]{Lemma}

\newtheorem{remark}[theorem]{Remark}

\numberwithin{equation}{section}

\DeclareMathOperator{\diam}{\textnormal{diam}}
\DeclareMathOperator{\f}{\textnormal{fill}}
\newcommand{\partialint}{\partial_{\textnormal{int}}}
\newcommand{\U}{\mathcal{U}}
\newcommand{\K}{\kappa}
\DeclareMathOperator{\Range}{\textnormal{Range}}

\begin{document}
\title{On the size of a finite vacant cluster of random interlacements with small intensity}
\author{Teixeira A.}
\maketitle

\begin{abstract}
In this paper we establish some properties of percolation for the vacant set of random interlacements, for $d \geqslant 5$ and small intensity $u$. The model of random interlacements was first introduced by A.S. Sznitman in \cite{sznitman}. It is known that, for small $u$, almost surely there is a unique infinite connected component in the vacant set left by the random interlacements at level $u$, see \cite{vladas} and \cite{teixeira}. We estimate here the distribution of the diameter and the volume of the vacant component at level $u$ containing the origin, given that it is finite. This comes as a by-product of our main theorem, which proves a stretched exponential bound on the probability that the interlacement set separates two macroscopic connected sets in a large cube. As another application, we show that with high probability, the unique infinite connected component of the vacant set is ``ubiquitous'' in large neighborhoods of the origin.
\end{abstract}

\section{Introduction}
\label{sec:intro}

In this paper we proceed with the study of the random interlacements introduced by A.S. Sznitman in \cite{sznitman}. This model is for instance related to the trace left by a random walk on the discrete torus $(\mathbb{Z}/N\mathbb{Z})^d$ ($d \geqslant 3$) and on the discrete cylinder $(\mathbb{Z}/N\mathbb{Z})^d \times \mathbb{Z}$ ($d \geqslant 2$) when the walk runs for times of order $N^d$ and $N^{2d}$ respectively, see \cite{benjamini} and \cite{sznitman_cil0}. Intuitively, random interlacements describe the microscopic `texture in the bulk' left by the random walk in these contexts, see \cite{david2} and \cite{sznitman_cil1}. In \cite{sznitman_cil3}, this model is the main ingredient to improve the upper bound on the disconnection time of a large discrete cylinder, and in \cite{sznitman_cil4} they are used to extend the the lower bound obtained in \cite{sznitman_cil2} to the case $d \geqslant 2$.

Loosely speaking, the interlacement at level $u$ (denoted by $\mathcal{I}^u$) is given by the trace left by a Poisson cloud of doubly infinite random walk trajectories in $\mathbb{Z}^d$, where $u$ controls the density of the cloud. The so-called vacant set at level $u$ (denoted with $\mathcal{V}^u$) is the complement of the interlacement, or in other words, the set of sites in $\mathbb{Z}^d$ which are not visited by any trajectory in this cloud. The random sets $\mathcal{I}^u$ are constructed simultaneously for all values of $u$ on the same probability space $(\Omega, \mathcal{A}, \mathbb{P})$.

Although we postpone the precise description of the process to Section~\ref{sec:review}, we state here a characterization of the law $Q^u$ of the indicator function of $\mathcal{V}^u$, regarded as a random element of $\{0,1\}^{\mathbb{Z}^d}$. Namely, $Q^u$ is the only probability measure on $\{0,1\}^{\mathbb{Z}^d}$ such that
\begin{equation}
\label{eq:Qucharacterized}
Q^u[Y_x = 1, \text{ for all } x \in K] = \exp \{-u \, \text{cap}(K)\}, \text{ for all finite sets } K \subset \mathbb{Z}^d,
\end{equation}
where $\text{cap}(K)$ denotes the capacity of $K$ (see (\ref{eq:cap})) and $(Y_x)_{x \in \mathbb{Z}^d}$ stand for the canonical coordinates on $\{0,1\}^{\mathbb{Z}^d}$, see Remark 2.2 2) of \cite{sznitman}.

Percolation of the vacant set of random interlacements presents a phase transition in the parameter $u$. More precisely, it is known that there is a $u_* \in (0,\infty)$, such that when $u < u_*$, $\mathcal{V}^u$ contains $\mathbb{P}$-a.s. an infinite connected component, see \cite{vladas} Theorem~3.4, and when $u > u_*$, $\mathcal{V}^u$ almost surely consists of finite clusters, see \cite{sznitman} Theorem~3.5. Moreover it is known that if an infinite connected component of the vacant set exists, it is almost surely unique, see \cite{teixeira} Theorem~1.1.

In this article we further investigate this model in the regime of small intensity $u$, and give a partial answer to the question posed in \cite{teixeira}, Remark 3.5 2). Denoting by $B(0,r)$ the closed ball with respect to the $l^\infty$ norm on $\mathbb{Z}^d$ with radius $r > 0$ and centered at the origin, our main Theorem~\ref{th:main} states that
\begin{equation}
\label{eq:maintheorem}
\begin{split}
& \text{for $d \geqslant 5$, there are $\bar u, \alpha >0$ such that, given $\gamma \in (0,1)$,}\\
& \mathbb{P} \left[ 
\begin{array}{c}
\text{for every pair of connected subsets of $B(0,N)$}\\
\text{with diameter at least $\gamma N$, there exists a path in}\\
\text{$\mathcal{V}^u \cap B(0,N+\gamma N)$ joining their boundaries}
\end{array} \right] \geqslant 1- c_1 e^{-c_2 N^\alpha},
\end{split}
\end{equation}
when $N \geqslant 0$ and $u \leqslant \bar u$. Here, $c_1$ and $c_2$ are positive constants solely depending on $\gamma$ and $d$, see Theorem~\ref{th:main}. 

\vspace*{4mm}

We will now underline the importance of this main result by stating some of its consequences, see also Section~\ref{sec:applications}. As a first application of (\ref{eq:maintheorem}), one can show the ubiquity of $\mathcal{C}_\infty^u$, the unique infinite component in $\mathcal{V}^u$.
More precisely, for $\gamma \in (0,1)$, $N \geqslant 1$ and $u \leqslant \bar u$, with the same notation as in (\ref{eq:maintheorem}), we show in Theorem~\ref{th:Cinftydense} that
\begin{equation}
\label{eq:stateCinftydense}
\begin{split}
& \mathbb{P}
\left[
\begin{array}{c}
\text{$\mathcal{C}_\infty^u$ intersects the boundary of every connected set}\\
\text{$C \subset B(0,N)$ with diameter at least $\gamma N$}
\end{array}
\right] \geqslant 1 - c_3 \exp(-c_4 N^\alpha),
\end{split}
\end{equation} 
where $c_3$ and $c_4$ are positive constants solely depending on $\gamma$ and $d$.

If $\mathcal{C}^u_0$ stands for the connected component of $\mathcal{V}^u$ containing the origin, with (\ref{eq:stateCinftydense}), we are able to control the diameter and the volume of $\mathcal{C}^u_0$, when this cluster is finite. This answers in part a question of \cite{sznitman}, Remark 4.4 3). Indeed in Theorems~\ref{th:taildiam} and \ref{th:tailvol}, we show that, for $d \geqslant 5$ and $u \leqslant \bar u$,
\begin{gather}
\label{eq:diamestimate}
\exp(-c_5(u) N) \leqslant \mathbb{P}\Big[\diam(\mathcal{C}_0^u) \geqslant N, |\mathcal{C}_0^u| < \infty \Big] \leqslant c_6 \exp(-c_7 N^\alpha) \,\,\, (N \geqslant 1),\\
\label{eq:volestimate}
c_8(u) \exp \big(-c_9V^{\frac{d-2}{d}} \log V \big) \leqslant \mathbb{P}\Big[ V \leqslant |\mathcal{C}_0^u| < \infty \Big] \leqslant c_6 \exp(-c_7 V^{\alpha/d}) \,\,\, (V \geqslant 1),
\end{gather}
where all the constants appearing above are positive and, except for $c_5$ and $c_8$, only depend on $d$. The estimate (\ref{eq:volestimate}) has the flavor of the open problem posed in \cite{benjamini}, Remark~4.7~1).

As another application of (\ref{eq:maintheorem}) we prove in Theorem~\ref{th:logcomp} that for $d \geqslant 5$, $u \leqslant \bar u$ and any $\epsilon > 0$ and integer $K \geqslant 1$,
\begin{equation}
\label{eq:statelogcomp}
\begin{split}
& \lim_{N \rightarrow \infty} N^K \mathbb{P} \left[
\begin{array}{c}
\text{some connected set of $\mathcal{V}^u \cap B(0,N)$}\\
\text{with diameter at least $(\log N)^{(1 + \epsilon)/\alpha}$ does not meet $\mathcal{C}^u_\infty$}
\end{array}
\right] = 0.
\end{split}
\end{equation}

\vspace*{4mm}

Let us mention that in the case of Bernoulli independent site percolation, similar results are already known to hold under weaker hypotheses. However, the techniques used in the Bernoulli context are not directly applicable to random interlacements, see Remark~\ref{rem:Bernoulli} 1) and 2). Some difficulties we mention here are the high dependence featured by the measure $Q^u$, see \cite{sznitman} (1.68), and the fact that for every $u > 0$ the interlacement at level $u$ is almost surely an infinite connected subset of $\mathbb{Z}^d$, see Corollary~2.3 of \cite{sznitman}. This property motivates the precise formulation of (\ref{eq:maintheorem}). As we further explain in Remark~\ref{rem:Bernoulli} 3), one can hope to find a vacant path joining the boundary of two large connected sets in $B(0,N)$, but not necessarily a vacant path joining the sets themselves.

\vspace*{4mm}

We now describe the strategy adopted to prove (\ref{eq:maintheorem}). In essence, the proof is based on the following two basic ingredients:
\begin{itemize}
\item [i)] if the interlacement separates two macroscopic components of a box, then in
many sub-boxes it also separates some macroscopic components. In other words,
the property of separating macroscopic components `cascades to finer scales',
\item [ii)] a fixed number of random walk paths can hardly separate macroscopic
components in a large box.
\end{itemize}
These claims are made precise and proved in Sections~\ref{sec:coarse} and \ref{sec:local}, respectively. In Section~\ref{sec:proof} we prove (\ref{eq:maintheorem}) using these results.

\vspace*{4mm}

Consider the following sequence of scales
\begin{equation}
L_\K = L_0(80L)^\K, \text{ for } \K \geqslant 0,
\end{equation}
where $L \geqslant 40, L_0 \geqslant 1$ are integers.

%

We now provide a short overview of the proof of \eqref{eq:maintheorem} and give an idea of the role plaid by the parameters $L$, $L_0$ and $u$.

Our aim is to bound the probability of the so-called separation event (in essence the complement of the event appearing in \eqref{eq:maintheorem}). In the above mentioned ingredient $i)$ of the proof, we show that separation `cascades down to finer scales'. This allows us to bound the probability of the separation event in a box at scale $\K$ (having diameter $L_\K$) by the probability that such separation occurs simultaneously in $2^\K$ well-spaced boxes at the bottom scale (each with diameter $L_0$). The bound on the latter probability has to be good enough to offset the number of possible choices for the boxes at the bottom scale. The combinatorial complexity of this choice roughly amounts to choosing a binary sub-tree of depth $\kappa$ in a rooted tree having $(\text{const} \cdot L)^{4d}$ descendants at each generation.

We thus need to control the probability of the simultaneous occurrence of separation events in each of the $2^\K$ boxes at the bottom scale $L_0$. For any such given collection of boxes, we first bound the mutual dependence of the separation events in each of them. For this purpose, we keep track of the number of excursions that the random walk trajectories (composing the random interlacements) perform between these boxes. We now choose a large enough $L$, consequently increasing the mutual distance between the boxes in this collection. In this fashion we are able to make the large deviation cost of observing too many excursions offset the combinatorial complexity of the choices of the $2^\K$ boxes. This step is delicate because increasing $L$ also increases this combinatorial complexity. For this competition to work in our favor, we need to impose the boxes at the bottom scale to receive an average number of excursions (say $a$) such that $a(d-2) > 4d$. We take $a = 100$ since this will do the job, see also Remark~\ref{rem:G100} 1). The dependence control described in this paragraph works for every $L_0$, once $u$ is chosen small enough depending on $L$ and $L_0$.

The previous step enables us to treat the separation event in the $2^\K$ boxes as roughly independent. We thus choose $L_0$ large enough so that the probability of these essentially independent $2^\K$ separation events in boxes at scale $L_0$ also offsets the combinatorial choice of the boxes. Now the above ingredient $ii)$ of the proof comes into play. When $L_0$ is chosen large enough (depending on $L$), a fixed number (as we said, $100$ will do the job) of independent random walk excursions can hardly separate components in a given box of size $L_0$. Now that $L_0$ is fixed, we can choose $u = u(L,L_0)$ small to make sure that observing an average of $100$ random walks per box in the collection above is indeed a large deviation as we described.

\bigskip

We now give a more precise description of the proof. For each depth $\K$, we partition $\mathbb{Z}^d$ into boxes of diameter $L_\K$ and label the boxes in this partition with a set of indices $I_\K$. For a given box (say indexed by $m \in I_\K$), we consider the random variable $\chi_m(\mathcal{I}^u)$, which loosely speaking indicates the separation of two macroscopic connected sets of this box by $\mathcal{I}^u$. We refer to (\ref{eq:chi}) for the precise definition.

As a reduction step, we prove that it is enough to establish (\ref{eq:maintheorem}) in the case where $\gamma = 2/3$ and $\gamma N$ is taken along the sequence $L_\K$, $\K \geqslant 0$. In other words, according to Proposition~\ref{prop:reduct}, in order to prove (\ref{eq:maintheorem}) we only need to obtain a bound (decaying exponentially in $2^\K$) on the probability that $\chi_m(\mathcal{I}^u) = 1$ when $m \in I_\K$.

We also rely on the concept of a skeleton, which captures the possible ways in which the separation event can propagate to finer scales. Roughly speaking a skeleton is a set $M$ of indices in the finest scale ($M \subset I_0$) satisfying some conditions on the distance between the boxes indexed by $M$, see Definition~\ref{def:skeleton}. The notion of skeleton resembles the Wiener criterion, see for instance \cite{lawler}, Theorem~2.2.5 p.55. The main purpose of this definition appears in (\ref{eq:boundexcursion}), where we derive a bound on the probability that a random walk, starting in one of the boxes of a skeleton, hits another box of the skeleton before escaping to infinity.

\vspace*{4mm}

The ingredient $i)$ of the proof, which we call `coarse graining' argument, is the content of Theorem~\ref{th:coarse2}, Section~\ref{sec:coarse}. Loosely speaking, this theorem states that
\begin{equation}
\label{eq:coarse}
\begin{array}{c}
\text{if $m \in I_\K$, and $\chi_m(\mathcal{I}^u) = 1$, there exists a skeleton $M \subset I_0$ with $\#M = 2^\K$} \\
\text{such that $\chi_{m'}(\mathcal{I}^u) = 1$ for all the indices $m' \in M$. Moreover,} \\
\text{the number of choices for such a skeleton is bounded by $((5 \cdot 80L)^{4d})^{2^\K}$}.
\end{array} 
\end{equation}
Hence, the problem is reduced to estimating the probability that the `separation event' ($\chi_{m'}(\mathcal{I}^u) = 1$) occurs simultaneously for the $2^\K$ indices $m'$ in a given skeleton $M$ as above. This bound has to be able to offset the combinatorial complexity factor $((5 \cdot 80L)^{4d})^{2^\K}$.

\vspace*{4mm}

Let us now indicate how the above bound is related to the second ingredient, which we call `local estimates'. For this, fix a skeleton $M \subset I_0$ and a collection of $2^\K$ boxes associated to indices $m$ in $M$. Loosely speaking, we use a large deviation estimate to bound the total number of excursions performed between different boxes of this collection by all the interlacement trajectories, see (\ref{eq:boundfirst}). Then we condition each of these excursions on their return and departure points from different neighborhoods of each box. This procedure will reduce our problem to the analysis of what happens in the surroundings of one fixed box of diameter $L_0$.

\vspace*{4mm}

The ingredient $ii)$ of the proof is obtained in Section~\ref{sec:local}. It can be summarized as follows:
\begin{equation}
\begin{array}{c}
\label{eq:statelocal}
\text{with high probability as $L_0$ grows, a fixed number} \\
\text{of independent random walk excursions do not produce a} \\
\text{`separation event' in the vicinity of a box of diameter $L_0$.} \\
\end{array}
\end{equation}
Moreover, this estimate is uniform on the points in which we condition these random walks to enter and exit a large neighborhood of the box. This is the content of Theorem~\ref{th:local}, Section~\ref{sec:local} and is the last piece to establish (\ref{eq:maintheorem}). Theorem~\ref{th:local} is the only part of the proof of our main result in which we need the hypothesis $d \geqslant 5$.

\vspace*{4mm}

Finally, let us outline of the proof of (\ref{eq:statelocal}). First we introduce the definition of a cut-point for a double infinite trajectory, see (\ref{eq:cuttime}). Loosely speaking, we regard each random walk excursion as a finite set of `sausages' connected by cut-points. Given two connected subsets $A_1$ and $A_2$ of a box, we show that:
\begin{equation}
\begin{array}{c}
\text{if the diameters of both $A_1$ and $A_2$ are big (when compared with}\\
\text{the diameter of each sausage), then we can connect the boundaries of}\\
\text{$A_1$ and $A_2$ by a path that avoids the corresponding random walk excursion},
\end{array}
\end{equation}
see Corollary~\ref{lem:chifillX}. Roughly speaking, we construct this path by ``traveling along the boundaries of the sausages''. Finally we show that with high probability (as $L_0$ grows) the diameters of the `sausages' are small when compared with $L_0$ and the excursions performed in the box $([0,L_0) \cap \mathbb{Z})^d$ are mutually far apart, so that they can be treated separately, see Lemma~\ref{lem:farapart}.

\vspace*{4mm}

This article is organized as follows.

In Section~\ref{sec:review} we give a precise description of the random interlacements and state some results which are used throughout the article.

In Section~\ref{sec:applications}, our main Theorem~\ref{th:main} is stated. In addition, we derive several applications of Theorem~\ref{th:main}. In Theorems~\ref{th:Cinftydense}, \ref{th:taildiam}, \ref{th:tailvol} and \ref{th:logcomp} we prove (\ref{eq:stateCinftydense}), (\ref{eq:diamestimate}), (\ref{eq:volestimate}) and (\ref{eq:statelogcomp}) respectively.

In Section~\ref{sec:proof} we prove the main Theorem~\ref{th:main} assuming Theorems~\ref{th:coarse2} and \ref{th:local} and Lemma~\ref{lem:indicespath}, which are proved in the subsequent sections.

The main result of Section~\ref{sec:coarse} is Theorem~\ref{th:coarse2}, which implements the `coarse graining' argument (\ref{eq:coarse}) used to reduce the problem to a microscopic scale.

In Section~\ref{sec:local} we describe the local picture of the process (see (\ref{eq:statelocal})). This is the content of Theorem~\ref{th:local}.

\vspace*{4mm}

Finally we comment on our use of constants. Throughout this article, $c$ and $c'$ will be used to denote positive constants depending only on $d$ (except when explicitly mentioned), which can change from place to place. We write $c_1, c_2, \dots$ for fixed positive constants (also depending only on $d$), which refer to their first appearance in the text.

{\bf Acknowledgments} - We are grateful to Alain-Sol Sznitman for important suggestions and
encouragement.

\section{A brief review of random interlacements}
\label{sec:review}

In this section we introduce some notation and describe the model of random interlacements. In addition, we recall some useful facts concerning the model.

For $a \in \mathbb{R}$, we write $\lfloor a \rfloor$ for the largest integer smaller or equal to $a$ and recall that
\begin{equation}
\label{eq:floorconvex}
\lfloor t a+(1-t)b \rfloor \in [\min \{a,b\},\max \{a,b\}], \text{ for all } a,b \in \mathbb{Z} \text{ and } t \in [0,1].
\end{equation}

We denote by $\{{\pmb e}_j\}_{j = 1, \dots, d}$ the canonical basis of $\mathbb{R}^d$ and write $\{\pi_j\}_{j = 1, \dots, d}$ for the corresponding orthogonal projections. For $y \in \mathbb{R}^d$ we denote by $\text{floor}(y)$ the element $x$ of $\mathbb{Z}^d$ such that $\lfloor \pi_j(y) \rfloor = \pi_j (x)$ for $j = 1,\dots, d$. Given $x, y \in \mathbb{R}^d$ we write $x \perp y$ if they are orthogonal for the usual scalar product.

We let $\lVert \cdot \rVert_\infty$ and $\lVert \cdot \rVert$ respectively denote the $l^\infty$ and the $l^1$ norms on $\mathbb{R}^d$ and $B(x,r)$ stand for the closed $l^\infty$-ball in $\mathbb{Z}^d$, i.e. $\{y \in \mathbb{Z}^d; \lVert x-y \rVert_\infty \leqslant r\}$. We say that two points $x,y \in \mathbb{Z}^d$ are neighbors if $\lVert x-y \rVert = 1$ (we also write $x \leftrightarrow y$) and if $\lVert x - y \rVert_\infty = 1$ we say that $x$ and $y$ are $*$-neighbors (and write $x \overset{*}{\leftrightarrow} y$). These definitions respectively induce the notions of connectedness and $*$-connectedness in $\mathbb{Z}^d$.

If $K \subset \mathbb{Z}^d$, we denote by $K^c$ its complement, by $|K|$ its cardinality and by $B(K,r)$ the $r$-neighborhood of $K$ for the $l^\infty$-distance, i.e. the union of the balls $B(x,r)$ for $x \in K$. The diameter of $K$ (denoted by $\diam(K)$) is the supremum of $\lVert x - y \rVert_\infty$ with $x,y \in K$. We define the boundary $\partial K$ (respectively the $*$-boundary $\partial^* K$) by $\{x \in K^c; x \leftrightarrow y \text{ for some } y \in K\}$ (respectively by $\{x \in K^c; x \overset{*}{\leftrightarrow} y \text{ for some } y \in K\}$). Analogously, we define the interior boundary $\partialint K = \{x \in K; x \leftrightarrow y \text{ for some } y \in K^c\}$. We also write $\overline K = K \cup \partial K$ and if $K$ is finite, we denote by $\text{sbox}(K)$ the smallest box containing $K$ (by box we mean a set of type $[a_1,b_1] \times \dots \times [a_d,b_d] \cap \mathbb{Z}^d$). For $K, K' \subset \mathbb{Z}^d$, the distance $d(K,K')$ is given by $\inf \{\lVert x-y \rVert; x \in K, y \in K'\}$ while $d_\infty(K,K') = \inf\{\lVert x - y \rVert_\infty; x \in K, y \in K'\}$ .

We will need the following
\begin{definition}
\label{def:fill}
Given a finite set $A \subset \mathbb{Z}^d$, we define $\f(A)$ as the complement of the unique unbounded connected component of $A^c$.
\end{definition}

Two important features of the set $\f(A)$ are stated in the following theorems.
\begin{gather}
\label{eq:stfillconn}
\text{If $A$ is connected, $\partial^* \f(A)$ is connected.}\\
\label{eq:fillstconn}
\text{If $A$ is connected, $\partial \f(A)$ is $\ast$-connected.}
\end{gather}
For the proofs, see for instance, \cite{kesten_0}, Lemma (2.23) p.139 and \cite{pisztora}, Lemma 2.1.

During this article the term \textit{path} (respectively $*$-\textit{path}) always denote finite, nearest neighbor (resp. $\ast$-neighbor) paths, i.e. some $\tau:\{0,\dots,n\} \rightarrow \mathbb{Z}^d$ such that $\tau(l) \leftrightarrow \tau(l+1)$ (resp. $\tau(l) \overset{*}{\leftrightarrow} \tau(l+1)$) for $l = 0,\dots, n-1$. In this case we say that the length of $\tau$ is $n$ and denote it by $N_\tau$.

We denote with $W_+$ and $W$ the spaces of infinite, respectively doubly infinite, transient trajectories
\begin{equation}
\begin{split}
W_+ = \Big\{w:\mathbb{Z}_+ \rightarrow \mathbb{Z}^d; w(l) \leftrightarrow w(l+1), \text{ for each } l \geqslant 0 \text{ and } \lVert w(l) \rVert \xrightarrow[l \rightarrow \infty]{} \infty \Big\},\\
W = \Big\{w:\mathbb{Z} \rightarrow \mathbb{Z}^d; w(l) \leftrightarrow w(l+1), \text{ for each } l \in \mathbb{Z} \text{ and } \lVert w(l) \rVert \xrightarrow[|l| \rightarrow \infty]{} \infty \Big\}.
\end{split}
\end{equation}
We endow them with the $\sigma$-algebras $\mathcal{W}_+$ and $\mathcal{W}$ generated by the coordinate maps $\{X_n\}_{n \in \mathbb{Z}_+}$ and $\{X_n\}_{n \in \mathbb{Z}}$. For $w \in W_+$ (or $W$), we write $X_{[a,b]}$ for the set $\{X_n; n \in [a,b]\}$ and analogously for $X_{(a,b]}$, $X_{[a,b)}$ and $X_{(a,b)}$.

We state a useful link between $\f(A)$ and paths in $W_+$. Namely,
\begin{lemma}
\label{eq:charfill}
\begin{equation*}
\f(A) = \Big\{z \in \mathbb{Z}^d; \text{ for all $w_+ \in W_+$ with $X_0(w_+) = z$, $\Range(w_+) \cap A \neq \varnothing$} \Big\}.
\end{equation*}
\end{lemma}
\textit{Proof.} We first show that the complement of the above set is included in $(\f(A))^c$. Indeed, the existence of $w_+ \in W_+$ with $\Range(w_+) \cap A = \varnothing$ implies that $z = X_0(w_+)$ belongs to the infinite connected component of $A^c$ and hence to $\f(A)^c$. Conversely, since $\f(A)^c$ is infinite and connected, given $z \in \f(A)^c$ one can find a trajectory $w_+ \in W_+$ starting in $z$ which is disjoint from $\f(A)$. The Lemma~\ref{eq:charfill} follows. $\qquad \square$

We also introduce the entrance and exit times of a finite set $K \subset \mathbb{Z}^d$
\begin{equation}
\label{eq:entranceexit}
\begin{array}{c}
H_K(w) = \inf \{k \in \mathbb{Z}_{(+)}; X_k(w) \in K\}, \text{ for } w \in W_{(+)}, \\
T_K(w) = \inf \{k \geqslant 0; X_k(w) \in K^c \}, \text{ for } w \in W \text{ or } W_+,
\end{array}
\end{equation}
and for $w \in W_+$, we define the hitting time of $K$
\begin{equation}
\label{eq:hittingtime}
\widetilde{H}_K(w) = \inf \{k \geqslant 1; X_k(w) \in K\}.
\end{equation}

Let $\theta_k:W \rightarrow W$ stand for the time shift given by $\theta(w)(\cdot) = w(\cdot + k)$ (where $k$ could also be a random time). Given finite sets $\Sigma \subset \widetilde \Sigma \subset \mathbb{Z}^d$, we consider on $W_+$ and $W$ the sequence or returns to $\Sigma$ and departures from $\widetilde \Sigma$
\begin{equation}
\label{eq:departreturn}
\begin{aligned}
R_1 &= H_\Sigma & D_1 & = R_1 + T_{\widetilde \Sigma} \circ \theta_{R_1}\\
R_n &= D_{n-1} + H_\Sigma \circ \theta_{D_{n-1}} & D_n & = R_n + T_{\widetilde \Sigma} \circ \theta_{R_n} \dots
\end{aligned}
\end{equation}
Notice that the stopping time $T_K$ is also defined in (\ref{eq:entranceexit}) for trajectories in $W$.

For $x \in \mathbb{Z}^d$, (recall that $d \geqslant 3$) we can define the law $P_x$ of a simple random walk starting on $x$ on the space $(W_+,\mathcal{W}_+)$. If $\rho$ is a measure on $\mathbb{Z}^d$, we write $P_\rho = \sum_{x\in \mathbb{Z}^d} \rho(x) P_x$. Their expectations are respectively denoted by $E_x$ and $E_\rho$. In some calculations, we may consider different dimensions, in this case we will explicitly write $P^d_x$ to avoid confusion.


We need some estimates in the hitting probability of a given set. First, let us define the Green function
\begin{equation}
g(x,y) = \sum_{n\geqslant 0} P_x[X_n = y], \text{ for } x,y \in \mathbb{Z}^d.
\end{equation}
We refer to \cite{lawler}, Theorem~1.5.4 p.31 for the following estimate
\begin{equation}
c' \frac{1}{1+|x-y|^{d-2}} \leqslant g(x,y) \leqslant c \frac{1}{|x-y|^{d-2}}, \text{ for } x,y \in \mathbb{Z}^d.
\end{equation}
We will use the following inequalities:
\begin{equation}
\label{eq:boundshitting}
\sum_{y\in K} g(x,y) \Big/ \sup_{z\in K} \big(\sum_{y \in K} g(z,y)\big) \leqslant P_x[H_K < \infty] \leqslant \sum_{y\in K} g(x,y) \Big/ \inf_{z\in K} \big(\sum_{y \in K} g(z,y)\big),
\end{equation}
see \cite{sznitman} (1.9). They follow by considering the bounded martingale $\sum_{y\in K} g(X_{n \wedge H_K},y)$ and remarking that it converges in $L^1(P_x)$ towards $\mathds{1}_{\{H_K < \infty\}} \sum_{y \in K} g(X_{H_K},y)$. The equality between the starting value of this martingale and the expectation (with respect to $P_x$) of its limit leads to the inequalities above.

Using \eqref{eq:boundshitting} and Theorem~1.5.4 of \cite{lawler}, p.31, we conclude that, for $x$ such that $|x| \geqslant 2r$,
\begin{equation}
\label{eq:hitfarball}
P_x[H_{B(0,r)} < \infty] \leqslant \frac{c r^d / (|x|-r)^{d-2}}{\sum_{h=1}^{r/2} c'' h^{d-1}/h^{d-2}} \leqslant c \Big(\frac{r}{|x|} \Big)^{d-2}
\end{equation}
Note that in the first inequality, we have split the last sum appearing in \eqref{eq:boundshitting} according to the $l^\infty$ distance from $z$ to $y$. See also (2.16) and Proposition~2.2.2 in \cite{lawler}, p.53.

Moreover, using the invariance principle, we obtain that
\begin{equation}
\label{eq:hitball}
\text{for $r \geqslant 1$, if } \lVert x \rVert_\infty \geqslant \tfrac{3r}{2},\text{ then } P_x[H_{B(0,r)} < \infty] \leqslant c_{10} < 1.
\end{equation}

We introduce, for a finite $K \subset \mathbb{Z}^d$, the equilibrium measure
\begin{equation}
e_K(x) = \mathds{1}_{x \in K} P_x[\widetilde H_K = \infty], \text{ for } x \in \mathbb{Z}^d,
\end{equation}
the capacity of $K$
\begin{equation}
\label{eq:cap}
\text{cap}(K) = e_K(\mathbb{Z}^d)
\end{equation}
and the normalized equilibrium measure
\begin{equation}
\label{eq:normalizedeK}
\overline e_K (x) = e_K(x)/ \text{cap}(K), \text{ for } x \in \mathbb{Z}^d.
\end{equation}

We mention the following bound on the capacity of a ball of radius $r$
\begin{equation}
\label{eq:estimatecap}
\text{cap}(B(0,r)) \leqslant c r^{d-2}, \text{ see \cite{lawler} (2.16), p.53.}
\end{equation} 

Let $W^*$ stand for the space of doubly infinite trajectories in $W$ modulo time shift,
\begin{equation}
W^* = W/\sim \mbox{, where } w \sim w' \mbox{ if } w(\cdot) = w'(k + \cdot), \mbox{ for some } k \in \mathbb{Z},
\end{equation}
endowed with the $\sigma$-algebra
\begin{equation}
\mathcal{W}^* = \{ A \subset W^* (\pi^*)^{-1}(A) \in \mathcal{W}\},
\end{equation}
which is the largest $\sigma$-algebra making the canonical projection $\pi^*:W \rightarrow W^*$ measurable. For a finite set $K \subset \mathbb{Z}^d$, we denote as $W_K$ the set of trajectories in $W$ which meet the set $K$ and define $W^*_K = \pi^*(W_K)$.

Now we are able to describe the intensity measure of the Poisson point process which governs the random interlacements.

For a finite set $K \subset \mathbb{Z}^d$, we consider the measure $Q_K$ in $(W, \mathcal{W})$ supported in $W_K$ such that, given $A,B \in \mathcal{W}_+$ and $x \in K$,
\begin{equation}
\label{eq:QK}
Q_K[(X_{-n})_{n \geqslant 0} \in A, X_0 = x, (X_n)_{n \geqslant 0} \in B] = P_x[A,\tilde H_K = \infty]P_x[B].
\end{equation}
Theorem~1.1 of \cite{sznitman} establishes the existence of a unique $\sigma$-finite measure $\nu$ in $W^*$ such that,
\begin{equation}
\label{eq:nuQK}
1_{W_K^*}\cdot \nu = \pi^* \circ Q_K, \text{ for any finite set $K \subset \mathbb{Z}^d$}
\end{equation}

We then introduce the spaces of point measures on $W^* \times \mathbb{R}_+$ and $W_+ \times \mathbb{R}_+$
\begin{equation}
\label{eq:Omega}
\begin{split}
& \begin{split}
\Omega = \bigg\{\omega = \sum_{i\geqslant 1} \delta_{(w^*_i,u_i)}; & w^*_i \in W^*, u_i \in \mathbb{R}_+ \mbox{ and } \omega(W^*_K \times [0,u]) < \infty\\
&\mbox{ for every finite } K \subset \mathbb{Z}^d \mbox{ and } u\geqslant 0 \bigg\}.
\end{split}\\
&\begin{split}
M = \bigg\{\mu = \sum_{i \in I} \delta_{(w_i,u_i)}; & I \subset \mathbb{N}, w_i \in W_+, u_i \in \mathbb{R}_+ \mbox{ and }\\
&\omega(W_+ \times [0,u]) < \infty \mbox{ for every } u\geqslant 0 \bigg\},
\end{split}
\end{split}
\end{equation}
endowed with the $\sigma$-algebras $\mathcal{A}$ and $\mathcal{M}$ generated by the evaluation maps $\omega \mapsto \omega(D)$ for $D \in \mathcal{W}^* \otimes \mathcal{B}(\mathbb{R}_+)$
and $\mu \mapsto \mu(D)$ for $D \in \mathcal{W}_+ \otimes \mathcal{B}(\mathbb{R}_+)$. Here $\mathcal{B}(\cdot)$ denotes the Borel $\sigma$-algebra.

We let $\mathbb{P}$ be the law of a Poisson point process on $\Omega$ with intensity measure $\nu \otimes du$, where $du$ denotes the Lebesgue measure on $\mathbb{R}_+$. Given $\omega = \sum_i\delta_{(w^*_i,u_i)} \in \Omega$, we define \textit{the interlacement} and the \textit{vacant set} at level $u$ respectively as the random subsets of $\mathbb{Z}^d$:
\begin{gather}
\label{eq:Iu}
\mathcal{I}^u (\omega) = \bigg\{ \bigcup_{i; u_i \leqslant u} \textnormal{Range}(w^*_i) \bigg\} \mbox{ and}\\
\label{eq:Vu}
\mathcal{V}^u(\omega) = \mathbb{Z}^d \setminus \mathcal{I}^u(\omega).
\end{gather}
We introduce the critical value
\begin{equation}
\label{eq:u*}
u_* = \inf \{u \geqslant 0; \mathbb{P}[\mathcal{V}^u \text{ contains an infinite connected component}] = 0\},
\end{equation}
c.f. \cite{sznitman}, (0.13).

It is known that for all $d \geqslant 3$,
$$
0 < u_* < \infty,
$$
see \cite{sznitman}, Theorem~3.5 and \cite{vladas}, Theorem~3.4. Moreover, its is also proved that if existent, the infinite connected component of the vacant set must be unique, see \cite{teixeira}, Theorem~1.1.

For a finite set $K \subset \mathbb{Z}^d$, we define the law $\mathbb{P}_K$ on $(M,\mathcal{M})$ of a Poisson point process in $W_+ \times \mathbb{R}_+$ with intensity measure $P_{e_K} du$.

The point processes defined above are related by the following. Consider, for a finite set $K \subset \mathbb{Z}^d$, the map $s_K:W^*_K \rightarrow W_+$ defined as
\begin{equation}
\begin{array}{c}
s_K(w^*) \text{ is the trajectory starting where}\\
w^* \text{ enters $K$ and following $w^*$ step by step,}
\end{array}
\end{equation} 
as well as the map $\mu_K:\Omega \rightarrow M$ defined via
\begin{equation}
\begin{split}
\mu_K(\omega)(f) = \int_{W^*_K\times \mathbb{R}_+} f(s_K(w^*),u) \omega(dw^* du), \text{ for $\omega \in \Omega$},\\
\end{split}
\end{equation}
for $f$ non-negative, measurable in $W_+ \times \mathbb{R}_+$.

One can prove that
\begin{equation}
\label{eq:PKmuP}
\mathbb{P}_K \text{ is the law of $\mu_K$ under $\mathbb{P}$},
\end{equation}
see \cite{sznitman} Proposition 1.3. And defining the point process $\mu_{K,u}$ on $\Omega$ with values in the set of finite point measures on $(W_+,\mathcal{W}_+)$
\begin{equation}
\label{eq:muKu}
\mu_{K,u}(\omega)(dw) = \mu_K(\omega)(dw\times [0,u]), \text{ for } \omega \in \Omega,
\end{equation}
we have
\begin{equation}
\label{eq:expexpectation}
\begin{array}{c}
\mathbb{E}[\exp\{- \langle \mu_{K,u},g \rangle \} ]= \exp\{u E_{e_K}[e^{- g} -1]\}, \\
\text{ for every non-negative $\mathcal{W}_+$-measurable function $g$}.
\end{array}
\end{equation}
To see why this holds, define $f: W_+ \times \mathbb{R}_+ \to \mathbb{R}_+$ by $f(w, v) = g(w)1_{\{v \leqslant u\}}$ and use (1.20) and (1.43) of \cite{sznitman}.

\section{The main result and some applications}
\label{sec:applications}

In this section we state our main result, Theorem~\ref{th:main} which translates (\ref{eq:maintheorem}) and, although we postpone the proof of this theorem to the next section, we now establish some of its consequences.

In Theorems~\ref{th:Cinftydense} we prove the existence of a unique infinite component in $\mathcal{V}^u$, for $d \geqslant 5$ and $u$ small enough, and show that with high probability this component is `dense' in the sense of (\ref{eq:stateCinftydense}). Theorems~\ref{th:taildiam} and \ref{th:tailvol} respectively provide estimates on the tail distribution of the diameter and the volume of the vacant component containing the origin when it is finite, see (\ref{eq:diamestimate}) and (\ref{eq:volestimate}). Finally, we bound the probability of finding a vacant component contained in $B(0,N)$ with diameter at least $\log(N)^{2/\alpha}$, as in (\ref{eq:statelogcomp}). This appears in Theorem~\ref{th:logcomp}.

\begin{definition}
\label{def:separate}
We say that two subsets $A_1$ and $A_2$ of $\mathbb{Z}^d$ are separated by $\U$ in $B$ if
\begin{enumerate}
\item $d(A_1,A_2) > 1$ and
\item every path in $B$ joining $\partial A_1$ to $\partial A_2$ meets $\U$.
\end{enumerate}
\end{definition}

We remark that
\begin{equation}
\label{eq:monotoneinB}
\begin{array}{c}
\text{if $B \subset B'$ and $\U \subset \mathbb{Z}^d$ separates $A_1$ from $A_2$ in $B'$,} \\
\text{then $A_1$ and $A_2$ are also separated by $\U$ in $B$.}
\end{array}
\end{equation} 

We now state the main result of the present article.

\begin{theorem}
\label{th:main}
($d \geqslant 5$) There are $\bar u > 0$ and $\alpha > 0$ such that, for every $0 < \gamma < 1$,
\begin{equation}
\label{eq:mainbound}
\mathbb{P} \Bigg[
\begin{array}{c}
\text{there exist connected sets $A_1, A_2 \subset B(0,N)$ with diameters}\\
\text{at least $\gamma N$ which are separated by $\mathcal{I}^u$ in $B(0,(1+\gamma)N)$}
\end{array}
\Bigg] < c_1 \cdot \exp(-c_2 N^\alpha),
\end{equation}
for all $N \geqslant 1$ and $u \leqslant \bar u$, where $c_1 = c_1(d,\gamma) > 0$ and $c_2 = c_2(d) > 0$. Note that the event appearing above decreases with $\gamma$.
\end{theorem}

Throughout the rest of this section, $\bar u$ and $\alpha$ stand for the values appearing above. Before going into the proof of the main theorem, we illustrate some of its applications.

\vspace*{4mm}

As a first consequence of Theorem~\ref{th:main}, we show that, for $u \leqslant \bar u$, there is almost surely a unique infinite connected component $\mathcal{C}_\infty^u$ of $\mathcal{V}^u$ and with overwhelming probability $\mathcal{C}_\infty^u$ neighbors all the macroscopic connected subsets of $B(0,N)$. More precisely,

\begin{theorem}
\label{th:Cinftydense}
($d \geqslant 5$) With $\bar u$ and $\alpha$ as in Theorem~\ref{th:main}, for every $u \leqslant \bar u$, there is $\mathbb{P}$-a.s., a unique infinite connected component $\mathcal{C}_\infty^u$ of $\mathcal{V}^u$, and for every $0 < \gamma < 1$,
\begin{equation}
\label{eq:Cinftydense}
\begin{split}
\mathbb{P}
\Bigg[
\begin{array}{c}
\text{$d(\mathcal{C}_\infty^u, C) \leqslant 1$ for all connected sets}\\
\text{$C \subset B(0,N)$ with $\diam(C) \geqslant \gamma N$, }
\end{array}
\Bigg] \geqslant 1 - c_3 \cdot \exp(-c_4 N^\alpha),
\end{split}
\end{equation}
for all $N \geqslant 1$. Again, $c_3 = c_3(d,\gamma) > 0$ and $c_4 = c_4(d,\gamma) > 0$.
\end{theorem}

\begin{remark}
\label{rem:u*pos}
\textnormal{Theorem~\ref{th:Cinftydense} yields an alternative proof that $u_* > 0$ when $d \geqslant 5$ (in fact $u_* \geqslant \bar u > 0$). The positivity of $u_*$ appeared first in \cite{sznitman}, Theorem~4.3 for $d \geqslant 7$. Later the result was established for all $d \geqslant 3$, see Theorem~3.4 of \cite{vladas}.}
\end{remark}


\vspace*{4mm}

\textit{Proof of Theorem~\ref{th:Cinftydense}.} We can suppose $\gamma < 1/4$. Applying Theorem~\ref{th:main} with $\gamma' = \gamma/2$ to a sequence of boxes centered at the origin with radius $2^jN$ ($j \geqslant 1$), we have
\begin{equation}
\label{eq:boundj}
\begin{split}
\mathbb{P}
\Bigg[&
\begin{split}
\text{for some $j \geqslant 1$, there are connected } & \text{sets $A_1, A_2 \subset B(0,2^jN)$ with}\\
\text{$\diam(A_1), \diam(A_2) \geqslant \gamma 2^{j-1} N$, which ar} & \text{e separated by $\mathcal{I}^u$ in $B(0,2^{j+1}N)$}
\end{split}
\Bigg]\\
&\leqslant \sum_{j=1}^\infty c_1 \exp(-c_2 (2^j N)^\alpha) \leqslant c_3 \cdot \exp(-c_4 N^\alpha).
\end{split}
\end{equation}

Now we prove that for $N \geqslant 10/\gamma$,
\begin{equation}
\label{eq:complemevent}
\begin{array}{c}
\text{in the complement of the event appearing in (\ref{eq:boundj}),}\\
\text{every connected set $C \subset B(0,N)$ with $\diam(C) \geqslant \gamma N$,}\\
\text{is neighbor of an infinite connected component of $\mathcal{V}^u$}.
\end{array}
\end{equation}

Once we establish the statement above, the $\mathbb{P}$-a.s. existence of the infinite connected component will follow from \cite{sznitman}, (2.4). And its uniqueness will be a consequence of \cite{teixeira}, Theorem~3.1.

Fix a connected set $C \subset B(0,N)$ with $\diam(C) \geqslant \gamma N$ and suppose we are in the complement of the event in (\ref{eq:boundj}). Taking $j=1$, we conclude that $C$ is not separated from $\partialint B(0,2N)$ by $\mathcal{I}^u$ in $B(0,4N)$. This implies that we can find a path in $\mathcal{V}^u \cap B(0,2N)$ starting at $\partial C$ and ending in $\partialint B(0,2N-1)$ so that its diameter is at least $N/2 \geqslant 2 \gamma N$. Let $C_1$ denote the range of this path.

Now suppose we have constructed connected sets $C_1 \subset C_2 \subset \dots \subset C_k$ in $\mathcal{V}^u \cap B(0,2^kN)$ with diameters at least $2 \gamma N, \dots, 2^k\gamma N$ respectively. We take $j = k+1$ in the complement of the event in (\ref{eq:boundj}) to conclude that there is a path in $\mathcal{V}^u \cap B(0,2^{k+1}N)$ connecting $\partial C_k$ to $\partialint B(0, 2^{k+1}N - 1)$. Hence, defining the connected set $C_{k+1}$ as the union of $C_k$ and the range of this path, we have $C_k \subset C_{k+1} \subset B(0,2^{k+1}N)$ and $\diam(C_{k+1}) \geqslant 2^{k-1} N \geqslant 2^{k+1} \gamma N$.

Letting $\bar C$ be $\cup_{j \geqslant 1}C_j$, we obtain an infinite connected subset of $\mathcal{V}^u$, which intersects $\partial C$ (since it contains $C_1$). This proves (\ref{eq:complemevent}) implying Theorem~\ref{th:Cinftydense}. $\quad \square$

\vspace*{4mm}

As another application of Theorem~\ref{th:main}, we consider $\mathcal{C}_0^u$, the connected component of $\mathcal{V}^u$ containing the origin and bound the tail of $\diam(\mathcal{C}_0^u)$, when $\mathcal{C}_0^u$ is finite.
\begin{theorem}
\label{th:taildiam}
($d \geqslant 5$) Let $\bar u$ and $\alpha$ be as in Theorem~\ref{th:main}. For every $u \leqslant \bar u$,
\begin{equation}
\exp(-c_5(u) N) \leqslant \mathbb{P}\Big[\diam(\mathcal{C}_0^u) \geqslant N, |\mathcal{C}_0^u| < \infty \Big] \leqslant c_6 \cdot \exp(-c_7 N^\alpha),
\end{equation}
where all the constants but $c_5 = c_5(d,u) > 0$ depend only on $d$.
\end{theorem}

\textit{Proof.} The upper bound follows from Theorem~\ref{th:Cinftydense}.

To prove the lower bound, we estimate the probability that $\mathcal{C}_0^u$ is precisely the segment $I_N = \{ j{\pmb e}_1; 0 \leqslant j \leqslant N\}$, where ${\pmb e}_1$ is the first vector in the canonical basis of $\mathbb{R}^d$, see below (\ref{eq:floorconvex}).

Define $I = \{ j{\pmb e}_1; j \in \mathbb{Z}\}$. Using the transience of the $(d-1)$-dimensional simple random walk, one concludes that for $x \in \partial I$
\begin{equation}
\label{eq:INIv}
P_x[\widetilde H_{\overline{I}_N} = \infty ] \geqslant P_x[\widetilde H_{\overline{I}} = \infty] = c > 0.
\end{equation}
and straightforwardly that
\begin{equation}
\label{eq:segmentscape}
P_x[\widetilde H_{\overline{I}_N} = \infty] \geqslant c > 0, \text{ for all $x \in \partial I_N$}.
\end{equation}

For $x \in \partial I_N$, write $W^{*,x} = \{w^* \in W^*; \text{$w^*$ enters $\overline{I}_N$ at $x$ and then leaves $\overline{I}_N$ forever} \}$.
Using (\ref{eq:segmentscape}) (\ref{eq:QK}) and (\ref{eq:nuQK}), we conclude that $\nu(W^{*,x}) \geqslant c^2$, for every $x \in \partial I_N$.

Since $\{W^{*,x}\}_{x \in \partial I_N}$ and $W^*_{I_N}$ are pairwise disjoint, the Poisson point processes obtained by restricting $\omega$ to these sets are independent, so that
\begin{equation}
\begin{split}
\mathbb{P} & \Big[\diam(\mathcal{C}_0^u) \geqslant N, |\mathcal{C}_0^u| < \infty\Big] \geqslant \mathbb{P}\big[\mathcal{C}_0^u = I_N\big]\\
& \geqslant \mathbb{P}\Big[\omega(W^{*,x} \times [0,u]) \geqslant 1, \text{ for all $x \in \partial I_N$ and } \omega(W_{I_N}^* \times [0,u]) = 0\Big] \\
& \geqslant \left(1-\exp\big(-u \cdot c^2 \big) \right)^{2dN} \exp(-u \cdot \text{cap}(I_N))\\
& \overset{(\ref{eq:cap})}{\geqslant} \exp(-c(u))^{2dN} \exp(-u \cdot N) \geqslant \exp(-c_5(u) N).
\end{split}
\end{equation}
This concludes the proof of Theorem~\ref{th:taildiam}. $\quad \square$
\vspace*{4mm}

We now bound the tail of the distribution of $|\mathcal{C}_0^u|$ (the volume of the vacant cluster containing the origin) when this cluster is finite.

\begin{theorem}
\label{th:tailvol}
For $u \leqslant \bar u$ ($\bar u$ as in Theorem~\ref{th:main}),
\begin{equation}
c_8(u) \cdot \exp \big(-c_9V^{\frac{d-2}{d}} \log V \big) \leqslant \mathbb{P}\Big[ V \leqslant |\mathcal{C}_0^u| < \infty \Big] \leqslant c_6 \cdot \exp(-c_7 V^{\alpha/d}),
\end{equation}
with all the constants but $c_8 = c_8(d,u) > 0$ depending only on $d$.
\end{theorem}

\textit{Proof.} The upper bound follows from Theorem~\ref{th:taildiam} once we use the fact that for some $c > 0$, $\diam(A) \geqslant c |A|^{1/d}$.

We prove the lower bound by estimating the probability that the box $B = B(0,N)$ is contained in $\mathcal{V}^u$ while the sphere $\partialint B(0,5N)$ is contained in $\mathcal{I}^u$. This is a case where $N^d \leqslant |\mathcal{C}_0^u| < \infty$.

We quote \cite{sznitman}, Remark~2.5 2) for the following estimate
\begin{equation}
\label{eq:coverbox}
P_y \big[B(y,N) \subset X_{[0,T_{B(y,2N)}]} \big] \geqslant c \cdot \exp(-c' N^{d-2} \log N).
\end{equation}

For $x \in \partialint B(0,5N)$ (we assume without loss of generality that $\pi_1(x) = 5N$) we consider the projection of the random walk starting at $x$ in the first coordinate. Using a gamblers ruin argument, one sees that with probability at least $c/N$ one reaches $\partial B(0,10N)$ before $\widetilde H_{B(0,5N)}$. From (\ref{eq:hitball}) we obtain
\begin{equation}
\label{eq:escape5N}
\begin{split}
P_x[\widetilde H_{B(0,5N)} = \infty] \geqslant \frac{c}{N}.
\end{split}
\end{equation}

We want the set $\partialint B(0,5N)$ to be contained in $\mathcal{I}^u$ and this will be the case if $B(Nx,N)$ $\subset$ $\mathcal{I}^u$ for all $x$ with $\lVert x \rVert_\infty = 5$. We define
\begin{equation}
\begin{split}
W^{*,x} = \Big\{ w^* & \in W^*_{B(0,5N)}; \text{ $w^*$ hits $B(0,5N)$ in $Nx$, then it covers $B(Nx,N)$ before}\\
& \text{leaving $B(Nx,2N)$ and escape to infinity without meeting $B(0,N)$} \Big\},
\end{split}
\end{equation}
for $x$ such that $\lVert x \rVert_\infty = 5$.

Using (\ref{eq:QK}) and (\ref{eq:nuQK}), we estimate
\begin{equation}
\begin{split}
\label{eq:boundnuWx}
\nu(W^{*,x}) \geqslant Q_{B(0,5N)}\big[ & X_0 = Nx, H_{B(0,N)} = \infty, B(Nx,N) \subset X_{ [0,T_{B(Nx,2N)}]}\big]\\
\overset{(\ref{eq:hitball}), (\ref{eq:escape5N}) , (\ref{eq:coverbox})}{\geqslant} & \frac{c}{N} \cdot \exp(-c'N^{d-2} \log N) \geqslant c \cdot \exp(-c'N^{d-2} \log N).
\end{split}
\end{equation}

Finally, since the sets $\{W^{*,x}\}_{x;\lVert x \rVert_\infty=5}$ and $W_{B(0,N)}^*$ are pairwise disjoint,
\begin{equation}
\begin{split}
\mathbb{P} & [N^d \leqslant |\mathcal{C}_0^u| < \infty] \geqslant \mathbb{P}[B(0,N) \subset \mathcal{V}^u, \partialint B(0,5N) \subset \mathcal{I}^u]\\
& \geqslant \mathbb{P}\Bigg[\bigg( \bigcap_{x \in \mathbb{Z}^d; \lVert x \rVert_\infty = 5}\omega(W^{*,x} \times [0,u]) \geqslant 1 \bigg) \cap \omega(W_{B(0,N)}^* \times [0,u]) = 0 \Bigg]\\
&\overset{(\ref{eq:boundnuWx})}{\geqslant} \Big( 1- \exp \big(-u \cdot c \cdot e^{-c'N^{d-2}\log N}\big) \Big)^{9^d \cdot 2d} \cdot c \exp (-u \cdot \text{cap}(B(0,N)))\\
& \overset{(\ref{eq:estimatecap})}{\geqslant} \Big( c(u) \cdot e^{-c'N^{d-2}\log N} \Big)^{9^d \cdot 2d} \cdot c \exp (-u \cdot c' \cdot N^{d-2}) \geqslant c(u) \exp(c' N^{d-2} \log N).
\end{split}
\end{equation}
And the proof is finished if one takes $N$ to be $\lfloor V^{1/d} \rfloor$. $\quad \square$

\vspace*{4mm}

\begin{remark}
\label{rem:intorus}
$\,$

\textnormal{1) The upper bound in Theorem~\ref{th:tailvol} has a flavor of the open problem posed in Remark~4.7 1) of \cite{benjamini}.}

\textnormal{2) Although the lower bound on Theorem~\ref{th:tailvol} is not expected to be sharp, it has a slower decay than the corresponding upper bound in the case of Bernoulli independent percolation, see for instance \cite{grimmett}, (8.66) p.216, c.f. Remark~1.1 of \cite{vladas}.} $\square$
\end{remark}

\vspace*{4mm}

As a last application of Theorem~\ref{th:Cinftydense}, we prove
\begin{theorem}
\label{th:logcomp}
($d \geqslant 5$) Let $\bar u$ and $\alpha$ be as in Theorem~\ref{th:main}. For every $u \leqslant \bar u$, $\epsilon > 0$ and $K \geqslant 1$,
\begin{equation}
\label{eq:logbound}
\lim_{N \rightarrow \infty} N^K \mathbb{P} \bigg[
\begin{array}{c}
\text{some connected set of $\mathcal{V}^u \cap B(0,N)$ with}\\
\text{diameter at least $(\log N)^{(1 + \epsilon)/\alpha}$ does not meet $\mathcal{C}^u_\infty$}
\end{array}
\bigg] = 0.
\end{equation}
\end{theorem}

\textit{Proof.} For a fixed $x \in B(0,N)$, one has
\begin{equation*}
\begin{split}
\mathbb{P} & \bigg[
\begin{array}{c}
\text{there is a finite connected component of $\mathcal{V}^u$ touching}\\
\text{$B(x,(\log N)^{(1 + \epsilon)/\alpha})$ with diameter larger or equal to $(\log N)^{(1 + \epsilon)/\alpha}$}
\end{array} \bigg]\\
&\leqslant \mathbb{P} \bigg[
\begin{array}{c}
\text{there is a connected set $C \subset B(x,2(\log N)^{(1 + \epsilon)/\alpha})$ which is}\\
\text{not neighbor of $\mathcal{C}^u_\infty$ and has diameter larger or equal to $\lfloor (\log N)^{(1 + \epsilon)/\alpha} \rfloor$}
\end{array} \bigg]\\
&\overset{(\ref{eq:Cinftydense})}{\leqslant} c \exp(-c' (\log N)^{(1 + \epsilon)}).
\end{split}
\end{equation*}
Summing over the points $x \in B(0,N)$, we obtain
\begin{equation}
\begin{split}
\mathbb{P} & \bigg[
\begin{array}{c}
\text{there is a connected subset of $\mathcal{V}^u$ in $B(0,N)$ with diameter}\\
\text{larger or equal to $(\log N)^{(1 + \epsilon)/\alpha}$ which is disjoint from $\mathcal{C}^u_\infty$}
\end{array} \bigg]\\
& \leqslant c N^{-c \log N + d}, \text{ and the claim (\ref{eq:logbound}) follows}. \quad \square
\end{split}
\end{equation}

\vspace*{4mm}

\begin{remark}
\label{rem:Bernoulli} \textnormal{Let us compare the results of this section with what is known to hold in the case of Bernoulli independent site percolation (where to every site one independently assigns the value $1$ with probability $p$ and $0$ with probability $1-p$).}

\textnormal{1) In the context of Bernoulli percolation, a result with a similar flavor to Theorem~\ref{th:main} can be proved. It is valid for any $d \geqslant 2$, any $p > p_c$ and providing an exponential bound instead of a stretched exponential, see for instance \cite{grimmett}, Lemma~(7.89) p.186. However, its proof strongly rely on the independence of the state of distinct sites, in contrast with the high dependence featured by the interlacement model, c.f. \cite{sznitman}, (1.68).}

\textnormal{2) In the Bernoulli independent case, one can use a Peierls-type argument to show that, for $p$ sufficiently close to one,}
\begin{equation}
\label{eq:Bernoullibound}
\begin{array}{c}
\textnormal{the probability that some $*$-connected component of 0's in $B(0,2N)$}\\
\textnormal{has diameter greater or equal to $N$ decays exponentially with $N$.}
\end{array}
\end{equation}
\textnormal{Together with (\ref{eq:stfillconn}), this provides a simple proof of the fact that: for this choice of~$p$, the probability that the diameter of the cluster of $1$'s containing the origin equals $N$ decays exponentially in $N$.}

\textnormal{Again this argument fails in the case of random interlacements. Actually, it is proved in \cite{sznitman}, Corollary~2.3, that}
\begin{equation}
\label{eq:Iuconnected}
\begin{array}{c}
\textnormal{$\mathcal{I}^u$ is $\mathbb{P}$-a.s an infinite connected subset of $\mathbb{Z}^d$ for all $u > 0$.}
\end{array}
\end{equation}
\textnormal{And according to (\ref{eq:Qucharacterized}), $\mathcal{I}^u$ meets $B(0,N)$ with overwhelming probability. Hence (\ref{eq:Bernoullibound}) does not hold for any $u > 0$ under the measure $Q^u$. See also \cite{sznitman}, Remark~2.5~2).}

\textnormal{3) As mentioned in the introduction, Theorem~\ref{th:main} does not hold true if one replaces $\partial A_1$ and $\partial A_2$ by $A_1$ and $A_2$ in Definition~\ref{def:separate}. We now give a brief justification for this claim.}

\textnormal{From the remark above we conclude that, with overwhelming probability as $N$ goes to infinity, one can find a self avoiding path $\tau$ contained in $\mathcal{I}^u$ connecting $B(0,N)$ to $\partialint B(0,2N)$. In this case, it is possible to extract two connected subsets $A_1$ and $A_2$ of $\Range(\tau)$ which have diameter at least $N/4$ and are far from each other. Since $A_1$ and $A_2$ are contained in $\mathcal{I}^u$, there is no path in $\mathcal{V}^u$ which joins these two sets. In the best case, we can hope to find a path in $\mathcal{V}^u$ which connects $\partial A_1$ to $\partial A_2$, as in Definition~\ref{def:separate}, see also the proof of Lemma~\ref{lem:avoid}.}

\textnormal{In the context of Bernoulli site percolation, we do not expect to need such restriction in the definition of separation, for a theorem analogous to Theorem~\ref{th:main} to hold for $p$ close enough to one.}

\textnormal{4) As an alternative definition of separation, one could for instance require that the sets $A_1$ and $A_2$ are disjoint from $\U$. This definition would be closer in spirit to the one appearing in \cite{grimmett}, Lemma~(7.89) p.186. However, in order to prove Theorem~\ref{th:main}, we will need the separation event to have a monotonicity property, see \eqref{eq:chimonotone}.}

\textnormal{Finding a suitable definition of separation is an important issue if one wishes to prove a theorem analogous to Theorem~\ref{th:main} that holds for any $u < u_*$.}
$\qquad \square$
\end{remark}

\section{Proof of the main result}
\label{sec:proof}

In this section we establish Theorem~\ref{th:main} using the Lemma~\ref{lem:indicespath} and the Theorems~\ref{th:coarse2} and \ref{th:local}, which are going to be proved in Sections~\ref{sec:coarse} and \ref{sec:local}.

Since the proof will follow a renormalization scheme, lets introduce the basic notation for the scales. For integers $L \geqslant 40, L_0 \geqslant 1$, we write
\begin{equation}
\label{eq:scales}
L_\K = L_0(80L)^\K, \text{ for $\K \geqslant 0$ integer},
\end{equation}
and define the set of indices in the scale $\K$
\begin{equation}
I_\K = \{\K\} \times \mathbb{Z}^d.
\end{equation}

Given $m = (\K,i) \in I_\K$, for $\K \geqslant 0$, we consider the box
\begin{equation}
\label{eq:boxCm}
C_m = (iL_\K + [0,L_\K)^d) \cap \mathbb{Z}^d
\end{equation}
and its $l$-th neighborhood
\begin{equation}
\label{eq:boxCml}
C_m^l = \bigcup_{j \in \mathbb{Z}^d; \lVert j \rVert_\infty < l} C(\K,i+j).
\end{equation}
Note that with this notation $C_m^1 = C_m$.

We introduce a definition which has a flavor of (1.6) in \cite{david}. For $m \in I_\K$, and $\U \subset \mathbb{Z}^d$, write
\begin{equation}
\label{eq:chi}
\begin{split}
\chi_m(\U) = \mathds{1}{
\bigg\{
\begin{array}{c}
\text{there exist connected sets $A_1, A_2 \subset C_m^2$, both with diameters}\\
\text{at least $L_\K /2$ which are separated by $\U$ in $C_m^3$}
\end{array}
\bigg\},
}
\end{split}
\end{equation}
recall the Definition~\ref{def:separate}.

Consider the function $\psi$ from $\{0,1\}^{\mathbb{Z}^d}$ into the set of subsets of $\mathbb{Z}^d$ that, given $\eta \in \{0,1\}^{\mathbb{Z}^d}$, is defined through $\psi(\eta) = \{x \in \mathbb{Z}^d; \eta(x) = 1\}$ (in other words, $\psi$ returns the subset of $\mathbb{Z}^d$ where $\eta$ is one). Since $\chi_m(\U)$ depends only on $\U$ intersected with $C_m^3$, one sees that $\{\chi_m \circ \psi = 1\}$ depends only on $Y_x$ for $x \in C_m^3$, therefore it is a cylinder event in the product $\sigma$-algebra of $\{0,1\}^{\mathbb{Z}^d}$.

It follows from the Definition~\ref{def:separate} that
\begin{equation}
\label{eq:chimonotone}
\chi_m \circ \psi \text{ is a non-decreasing function on $\{0,1\}^{\mathbb{Z}^d}$}.
\end{equation}

First we prove the following reduction step.
\begin{proposition}
\label{prop:reduct}
To prove Theorem~\ref{th:main}, it is suffices to show that:
\begin{equation}
\label{eq:reduct}
\begin{array}{c}
\text{there exist $\bar u > 0$, $L \geqslant 40$ and $L_0 \geqslant 1$ such that}\\
\mathbb{P}[\chi_{(\K,0)}(\mathcal{I}^u) = 1 ] \leqslant c \cdot 2^{-2^\K}, \text{ for every $u \leqslant \bar u$ and $\kappa \geqslant 0$}.
\end{array}
\end{equation}
\end{proposition}

\textit{Proof.} Fix $u \leqslant \bar u$ as above and recall the definition of $\text{sbox}(\cdot)$ above Definition~\ref{def:fill}. As an intermediate step towards (\ref{eq:mainbound}), we first prove that (\ref{eq:reduct}) implies that for any $0 < \gamma < 1$ and $\K \geqslant 1$,
\begin{equation}
\label{eq:reductinterm}
\begin{split}
& \mathbb{P} \Bigg[
\begin{split}
\text{there } & \text{exist $A_1, A_2 \subset C_{(\K,0)}$ with $\diam(A_1), \diam(A_2) \geqslant \tfrac{\gamma}{4} L_{\K - 1}$}\\
\text{which are} & \text{ connected and separated by $\mathcal{I}^u$ in $B(\text{sbox}(A_1 \cup A_2), \tfrac{\gamma}{4} L_{\K-1})$}
\end{split}
\Bigg]\\ & \leqslant c(L,L_0,\gamma) \, 2^{-L_\K^\alpha}, \text{with $\alpha = \alpha(L) > 0$}.
\end{split}
\end{equation}

Note that the condition on the diameters of $A_1$ and $A_2$ above is less restrictive than the condition appearing in (\ref{eq:chi}). Hence, on the event appearing in (\ref{eq:reductinterm}) there is no guarantee that $\chi_m (\mathcal{I}^u) = 1$ for some $m \in I_\K$ so that we cannot apply (\ref{eq:reduct}) with this $\K$ to prove (\ref{eq:reductinterm}). Instead, we choose an appropriate $\K_o \geqslant 1$, find an $\bar m \in I_{\K - \K_o}$ such that $\chi_{\bar m} (\mathcal{I}^u) = 1$ and then use (\ref{eq:reduct}). More precisely, let $\K_o = \K_o(\gamma) \geqslant 1$ be such that
\begin{equation}
\label{eq:n0gamma}
(80L)^{\K_o-1} > \tfrac{40}{\gamma},
\end{equation}
and note that, if $\K > \K_o$,
\begin{equation}
\label{eq:kkosmall}
L_{\K - \K_o} = \frac{L_{\K-1}}{(80L)^{\K_o-1}} < \frac{\gamma}{40} L_{\K-1}.
\end{equation}

On the event appearing in (\ref{eq:reductinterm}), one can find $i_1, i_2 \in \mathbb{Z}^d$ such that the boxes $C_{(\K - \K_o, i_1)}$ and $C_{(\K - \K_o, i_2)}$ are contained in $C_{(\K,0)}$ and intersect $A_1$ and $A_2$ respectively. By (\ref{eq:kkosmall}), $C_{(\K - \K_o, i_1)}^3, C_{(\K - \K_o, i_1)}^3 \subset B(\text{sbox}(A_1 \cup A_2), \tfrac{\gamma}{4} L_{\K-1}) \subset C_{(\K,0)}^2$. We now join $i_1$ and $i_2$ by a path $\tau$ in $\text{sbox}(i_1,i_2)$ and note that all boxes $C^3_{(\K - \K_o, \tau(n))}$ ($1 \leqslant n \leqslant N_\tau$) are contained in $B(\text{sbox}(A_1 \cup A_2), \tfrac{\gamma}{4} L_{\K-1})$. By Lemma~\ref{lem:indicespath} there is an index $\bar m = (\K - \K_o,\bar i)$ such that $C_{\bar m}^3 \subset B(\text{sbox}(A_1 \cup A_2), \tfrac{\gamma}{4} L_{\K-1}) \subset C_{(\K,0)}^2$ and $\chi_{\bar m}(\mathcal{I}^u) = 1$. As a result we have the bound
\begin{equation}
\label{eq:reducetoK}
\begin{split}
& \mathbb{P} \Bigg[
\begin{split}
\text{there } & \text{exist $A_1, A_2 \subset C_{(\K,0)}$ with $\diam(A_1), \diam(A_2) \geqslant \tfrac{\gamma}{4} L_{\K - 1}$}\\
\text{which are } & \text{connected and separated by $\mathcal{I}^u$ in $B(\text{sbox}(A_1 \cup A_2), \tfrac{\gamma}{4} L_{\K-1})$}
\end{split}
\Bigg]\\
& \qquad \leqslant \mathbb{P} \Bigg[
\begin{split}
\text{there } & \text{is some } \bar m = (\K-\K_o, \bar i), \text{ such that }\\
& C_{\bar m}^3 \subset C_{(\K,0)}^2 \text{ and } \chi_{\bar m}(\mathcal{I}^u) = 1
\end{split}
\Bigg] \\
& \qquad \leqslant (3 \cdot (80L)^{\K_o})^d \cdot \mathbb{P}[\chi_{(\K - \K_o,0)}(\mathcal{I}^u) = 1] \overset{(\ref{eq:reduct})}{\leqslant} c(L,\gamma) 2^{-2^{\K-\K_o}}.
\end{split}
\end{equation}

To obtain (\ref{eq:reductinterm}) we now set
\begin{equation}
\alpha = \frac{\log 2}{2 \log(80L)},
\end{equation}
and note that
\begin{equation}
\frac{\log(2^{\K-\K_o})}{\log L_{\K}} \overset{(\ref{eq:scales})}{=} \frac{(\K - \K_o) \log 2}{\K \log(80L) + \log (L_0)} \xrightarrow[\K \rightarrow \infty]{} \frac{\log 2}{\log (80L)} > \alpha,
\end{equation}
so that $2^{-2^{\K-\K_o}} \leqslant 2^{-L_\K^\alpha}$, for $\K$ larger or equal to some $\K(L_0,\gamma)$. Possibly increasing $c(L,\gamma)$ to some $c(L,L_0,\gamma)$ we obtain (\ref{eq:reductinterm}) from (\ref{eq:reducetoK}).

Finally, for a given $N$ (which we assume for the moment to be larger or equal to $L_0$) we choose a depth $\K(N) \geqslant 1$ such that $L_{\K(N)}$ is comparable with $N$, more precisely,
\begin{equation}
\label{eq:nN}
L_{\K(N)} > 2N \geqslant L_{\K(N)-1}.
\end{equation}

The boxes which appear in (\ref{eq:mainbound}) are centered at the origin (unlike the box $C_{(\K,0)}$ in (\ref{eq:reductinterm})). However, note that $B(0,N)$ and $B(0,(1+\gamma)N)$ can be respectively mapped (under the same translation) to the boxes $B_1 = [0,2N]^d \cap \mathbb{Z}^d$ and $B_2 = [-\gamma N, 2N + \gamma N]^d \cap \mathbb{Z}^d$. Hence, using that the law $\mathbb{P}$ is invariant under translation, see \cite{sznitman} Proposition~1.3, we conclude that it suffices to prove (\ref{eq:mainbound}) with $B(0,N)$ and $B(0,(1+\gamma)N)$ respectively replaced by $B_1$ and $B_2$. Note, by (\ref{eq:nN}), that
\begin{equation}
\label{eq:B1B2in}
B_1 \subset C_{(\K(N),0)},\quad \tfrac{\gamma}{4} L_{\K(N) - 1} \leqslant \tfrac{\gamma}{2}N \quad \text{and} \quad B(\text{sbox}(A_1 \cup A_2), \tfrac{\gamma}{4} L_{\K(N) - 1}) \subset B_2.
\end{equation}

Hence,
\begin{equation*}
\begin{split}
& \mathbb{P} \Bigg[
\begin{split}
 \text{the} & \text{re exist $A_1, A_2 \subset B_1$ with $\diam(A_1), \diam(A_2) \geqslant \gamma N$}\\
 & \text{which are connected and separated by $\mathcal{I}^u$ in $B_2$}
\end{split}
\Bigg]\\
& \overset{(\ref{eq:B1B2in}), (\ref{eq:monotoneinB})}{\leqslant} \mathbb{P} \Bigg[
\begin{split}
\text{there} & \text{ exist $A_1, A_2 \subset C_{(\K(N),0)}$ with $\diam(A_1), \diam(A_2) \geqslant \tfrac{\gamma}{4} L_{\K(N) - 1}$}\\
\text{which } & \text{are connected and separated by $\mathcal{I}^u$ in $B(\text{sbox}(A_1 \cup A_2), \tfrac{\gamma}{4} L_{\K(N)-1})$}
\end{split}
\Bigg] \\
& \overset{(\ref{eq:reductinterm})}{\leqslant} c(L,L_0,\gamma) \, 2^{-L_{\K(N)}^\alpha} \overset{(\ref{eq:nN})}{\leqslant} c_1(L,L_0,\gamma) \cdot \exp(-c_2 N^\alpha).
\end{split}
\end{equation*}
By possibly increasing $c_1$ we obtain (\ref{eq:mainbound}) for all $N \geqslant 1$. $\quad \square$


\vspace*{4mm}

In view of the above proposition, it suffices to show (\ref{eq:reduct}) in order to establish Theorem~\ref{th:main}. This will be done using the two main ingredients mentioned in the introduction, see (\ref{eq:coarse}) and (\ref{eq:statelocal}). To make this rough description precise, we will need the following
\begin{definition}
\label{def:skeleton}
A set of indices of boxes in the finest scale $M \subset I_0$ is said to be a \textit{skeleton} if, for any $m_0 \in M$, we have
\begin{equation}
\label{eq:skeleton}
\begin{aligned}
i) & \,\, \# \{m \in M; L \, L_h < d_\infty(C_{m_0}, C_m) \leqslant L \, L_{h+1}\} \leqslant 2^{h+1} \text{ for all $h \geqslant 0$ and} \\
ii) & \,\, \text{there is no box $m \in M \setminus \{m_0\}$ satisfying $d_\infty(C_{m_0}, C_m) \leqslant L \, L_0 $}.
\end{aligned}
\end{equation}
\end{definition}
This is a type of Wiener criterion, see for instance \cite{lawler}, Theorem~2.2.5, p.55.

\vspace*{4mm}

\textit{Proof of Theorem~\ref{th:main}.} The first ingredient of the proof of (\ref{eq:reduct}) is stated in the following theorem, which will be a direct consequence of Theorem~\ref{th:coarse2} proved in Section~\ref{sec:coarse}.

\begin{theorem}
\label{th:coarsefirst}
For each $m' \in I_\K$ ($\K \geqslant 0$), there is a family of skeletons $\mathcal{M}_{m'}$, such that
\begin{equation}
\label{th:coarse}
\begin{aligned}
\quad \textit{i)}&\,\, \text{for every skeleton $M \in \mathcal{M}_{m'}$, $\# M = 2^\kappa$},\\
\quad \textit{ii)}&\,\, \text{$\#\mathcal{M}_{m'} \leqslant ((5 \cdot 80L)^{4d})^{2^\K}$,}\\
\quad \textit{iii)}&\,\, \text{if $\chi_{m'}(\U) = 1$, we can find a skeleton $M \in \mathcal{M}_{m'}$} \\
\quad \quad &\text{(which may depend on $\U$), such that $\chi_m(\U) = 1$ for all $m\in M$.}
\end{aligned}
\end{equation}
\end{theorem}

The above statement will reduce the proof of (\ref{eq:reduct}) to the derivation of an appropriate bound on $\mathbb{P}[\cap_{m \in M} \{ \chi_m(\mathcal{I}^u) = 1\}]$ uniformly over the skeleton $M$. Roughly speaking, to obtain such bound we will analyze the excursions of the interlacement trajectories between neighborhoods of the boxes of the skeleton.

Given a skeleton $M$, we consider (see (\ref{eq:boxCml}) for the notation),
\begin{equation}
\label{eq:Sigmas}
\Sigma \overset{\text{def}}{=} \bigcup_{m\in M} C_m^5 \subset \widetilde \Sigma \overset{\text{def}}{=} \bigcup_{m\in M} C_m^{L/4},
\end{equation}
recalling that in (\ref{eq:boxCml}) we did not require $l$ to be an integer. By (\ref{eq:skeleton}) $ii)$, we conclude that for two distinct $m,m' \in M$, $d(C^{L/4}_m,C^{L/4}_{m'}) \geqslant LL_0/2 > 10$. Hence, for $x$ in $\widetilde \Sigma \cup \partial \widetilde \Sigma$, we can write $m(x)$ for the unique index $m \in M$ such that $x \in C^{L/4}_m \cup \partial C^{L/4}_m$.

Define the successive times of return to $\Sigma$ and departure from $\widetilde \Sigma$, $R_i$ and $D_i$ as in (\ref{eq:departreturn}) and note that on $\{X_0 \in \Sigma\}$, $R_1 = 0$. For $w \in W$, we define the number of excursions performed by $w$, as $g_M(w) = \sum_{l \geqslant 1} \mathds{1}_{\{R_l < \infty\}}$ and write
\begin{equation}
\label{eq:twoterms}
\begin{split}
\mathbb{P}\bigg[\bigcap_{m\in M} \{ \chi_m(\mathcal{I}^u) = 1\}\bigg] &\leqslant \mathbb{P}\Big[\langle \mu_{\Sigma,u}, g_M \rangle > 100 \cdot 2^\K\Big]\\
&\quad + \mathbb{P}\bigg[\langle \mu_{\Sigma,u}, g_M \rangle \leqslant 100 \cdot 2^\K, \bigcap_{m\in M} \{ \chi_m(\mathcal{I}^u) = 1\}\bigg].
\end{split}
\end{equation}

\vspace*{4mm}

To estimate the first term in the right-hand side of the above equation, we will use an exponential Chebychev-type inequality. With (\ref{eq:expexpectation}), this will amount to bounding the exponential moments of $g_M$ under $P_x$. This will be performed by choosing $L$ sufficiently large, consequently reducing the probability that a random walk starting on $\partial \widetilde \Sigma$ hits $\Sigma$ before escaping to infinity.

More precisely, for $y \in \partial \widetilde \Sigma$,
\begin{equation}
\label{eq:boundexcursion}
P_y[R_1 < \infty] \overset{(\ref{eq:hitfarball}), \text{(\ref{eq:skeleton})}}{\leqslant} P_y[H_{C_{m(y)}^5} < \infty] + \sum_{h=0}^\infty 2^{h+1} \cdot c \cdot \left( \frac{L_0}{L \cdot L_h} \right)^{d-2} \leqslant \frac{c_{11}}{L^{d-2}}.
\end{equation}

Using the Strong Markov property at time $R_2$, we obtain, for $x \in \Sigma$,
\begin{equation}
\begin{split}
\phi(x) & \mathrel{\mathop:}= E_x[e^{\lambda g_M}] \leqslant e^\lambda (1 + E_x[1\{R_2 < \infty\} \phi(X_{R_2})])\\
& \overset{(\ref{eq:boundexcursion})}{\leqslant} e^\lambda \Big(1+ \frac{c_{11}}{L^{d-2}} \sup_{z \in \Sigma} \phi(z) \Big).
\end{split}
\end{equation}
This inequality implies that
\begin{equation}
\label{eq:boundphi}
\text{when } e^\lambda \left(\frac{c_{11}}{L^{d-2}} \right) \leqslant 1/2, \text{ then } \sup_{z \in \Sigma} \phi(z) \leqslant 2 \cdot e^\lambda.
\end{equation}

Now choose
\begin{equation}
\label{eq:lambdachoice}
\lambda = \tfrac{1}{100} \log (2^{10} (5 \cdot 80L)^{4d})
\end{equation}
so that, c.f. (\ref{eq:skeleton}) $i)$,
\begin{equation}
\label{eq:etolambda}
e^{100\lambda} = 2^{10} (5 \cdot 80L)^{4d}.
\end{equation}

As a result, when we take $L \geqslant c$,
\begin{equation}
\label{eq:Lchoice}
\begin{array}{c}
e^{100\lambda}\left( \frac{c_{11}}{L^{d-2}} \right)^{100} \Bigg( = \big( 2^{10} (5 \cdot 80)^{4d} c_{11}^{100} \big) \left(\frac{L^{4d}}{L^{100(d-2)}} \right) \Bigg) < 2^{-100}, \\
\text{and the conclusion of (\ref{eq:boundphi}) holds.}
\end{array}
\end{equation}

Then, using (\ref{eq:expexpectation}), an exponential Chebychev-type inequality and $\text{cap}(\Sigma) \leqslant c2^\K L_0^{d-2}$, we get the desired bound on the first term of (\ref{eq:twoterms}), namely:
\begin{equation}
\label{eq:boundfirst}
\mathbb{P}[\langle \mu_{\Sigma,u}, g_M \rangle > 100 \, 2^\K] \leqslant \exp(-100 \lambda \, 2^\K + u E_{e_\Sigma} [e^{\lambda g_M} -1] ) \overset{(\ref{eq:boundphi})}{\leqslant} \left[ e^{-100\lambda + c\cdot u \cdot L_0^{d-2} e^\lambda} \right]^{2^\K} \negthickspace.
\end{equation}

The term in the right hand side of the equation above will be controlled at the very end of the proof of Theorem~\ref{th:main} using (\ref{eq:etolambda}) and choosing $u$ small enough.

\vspace*{4mm}

We now bound the second term of (\ref{eq:twoterms}). For this we need to invoke Theorem~\ref{th:local} of Section~\ref{sec:local}, that was referred as the second ingredient of the proof of Theorem~\ref{th:main}, see (\ref{eq:statelocal}). We need to introduce some notation.

For $x \in \partialint C_m^5$, $y \in \partial C_m^{L/4}$, recall that $P_x$-a.s. $R_1 = 0$ and $D_1 = T_{\widetilde \Sigma}$. Define
\begin{equation}
\label{eq:Pxy}
P_{x,y} [\,\cdot\,] = P_x[\,\cdot\, |X_{D_1} = y].
\end{equation}
Given $\vec{x} \in \Delta_1 = \{(x^i)_{i \leqslant G}; x^i \in \partialint C_{(0,0)}^5 \}$ and $\vec{y} \in \Delta_2 = \{(y^i)_{i \leqslant G}; y^i \in \partial C_{(0,0)}^{L/4} \}$, we define
\begin{equation}
\label{eq:Pvecxvecy}
P_{\vec{x},\vec{y}} = \bigotimes_{i \leqslant G} P_{x^i,y^i},
\end{equation}
and denote with $(X^i_n)_{i \leqslant G, n \geqslant 0}$ the canonical coordinates in $W_+^{\times G}$. We introduce the stopping time $D = T_{C_{(0,0)}^{L/4}}$.

We now need the second ingredient of the proof which is Theorem~\ref{th:local}. For the reader's convenience we state here this theorem and reefer to Section~\ref{sec:local} for its proof.
\newtheorem*{local}{Theorem~\ref{th:local}}
\begin{local}
\label{th:localfirst}
For $d \geqslant 5$, given $\epsilon > 0$, $G \geqslant 1$ and $L \geqslant 40$, for large enough $L_0 \geqslant 1$
\begin{equation}
\label{th:reftolocal}
\begin{array}{c}
{\displaystyle \sup_{\mbox{\fontsize{8}{8} \selectfont $\substack{\vec{x} \in \Delta_1, \vec{y} \in \Delta_2}$}}} P_{\vec{x},\vec{y}} \left [\chi_{(0,0)} \bigg( \bigcup_{i=1}^G X^i_{[0,D]} \bigg) = 1 \right] < \epsilon.
\end{array}
\end{equation}
\end{local}

The above result concerns only what happens in one fixed box in the finer scale. Loosely speaking, to bound the second term in the right hand side of (\ref{eq:twoterms}), we are going to condition all the interlacement trajectories which intersect $\Sigma$ on their return and departure points from $\Sigma$ and $\widetilde \Sigma$ (i.e., $X_{R_i}$ and $X_{D_i}$). We first need some further notation.

It is known that $\mu_{\Sigma, u}$ has the same distribution as $\sum_{1\leqslant i \leqslant S} \delta_{w^i}$, where $S \geqslant 0$ is a Poisson variable with parameter $u \, \text{cap}(\Sigma)$ and $w^i$ ($1 \leqslant i \leqslant S$) are i.i.d., $P_{\bar e_\Sigma}$-distributed and independent of $S$, recall the notation in (\ref{eq:normalizedeK}) and (\ref{eq:muKu}).

Write $X^1_l,\dots,X^S_l$ ($l \in \mathbb{Z}_+$) for the canonical coordinates of the trajectories $w^1,\dots,w^S$. We consider the number of excursions performed by $w_i$, $G_i = g_M(w^i)$ ($1 \leqslant i \leqslant S$) and for $1 \leqslant j \leqslant G_i$, we denote the range of the $j$-th excursion by $E_j^i = X_{[R_j,D_j]}^i$.

Given an index $m$ in the skeleton $M$, we collect the indices of all excursions performed in $C^{L/4}_m$:
\begin{equation}
\label{eq:Phim}
\Phi_m = \{(i,j); 1 \leqslant i \leqslant S, \,\, 1 \leqslant j \leqslant G_i, \,\, X_{R_j}^i \in C_m^5 \}.
\end{equation}

Note that, by (\ref{th:coarse})~\textit{i)}, when $\langle \mu_{\Sigma,u}, g_M \rangle \leqslant 100 \cdot 2^\K$, the set of $m' \in M$ with $|\Phi_{m'}| > 200$ has cardinality at must $2^{\K - 1}$, and hence
\begin{equation}
\label{eq:countM}
\begin{array}{c}
\text{when $\langle \mu_{\Sigma,u}, g_M \rangle \leqslant 100 \cdot 2^\K$, one can find some subset $M'$ of $M$}\\
\text{such that $|M'| = 2^{\K-1}$ and $|\Phi_{m'}| \leqslant 200$ for all $m' \in M'$.}
\end{array}
\end{equation}

Therefore, we can bound the last term of (\ref{eq:twoterms}) as follows:
\begin{equation*}
\begin{split}
& \mathbb{P}\bigg[ \langle \mu_{\Sigma,u}, g_M \rangle \leqslant 100 \cdot 2^\K, \bigcap_{m\in M} \{ \chi_m(\mathcal{I}^u) = 1\}\bigg]\\
& \overset{(\ref{eq:countM})}{\leqslant} \sum_{\substack{M'\subset M; \\ |M'| = 2^{\K-1}}} \mathbb{P} \bigg[ \bigcap_{m'\in M'} \Big(\{ \chi_{m'}(\mathcal{I}^u) = 1\} \cap \{ |\Phi_{m'}| \leqslant 200 \} \Big) \bigg] \\
& \overset{(\ref{eq:PKmuP})}{\leqslant} \sum_{\substack{M'\subset M; \\ |M'| = 2^{\K-1}}} \mathbb{P}_K \bigg[ \bigcap_{m'\in M'} \Big( \{ \chi_{m'}(\cup_{i \leqslant n} \Range(X^i)) = 1\} \cap \{ |\Phi_{m'}| \leqslant 200 \} \Big) \bigg]\\
& \leqslant 2^{2^\K} \sup_{\substack{M'\subset M; \\ |M'| = 2^{\K-1}}} \sum_{n\geqslant 0} \mathbb{P}[S=n] P_{\bar e_\Sigma}^{\otimes n} \bigg[ \bigcap_{m'\in M'} \Big( \{ \chi_{m'}(\cup_{i \leqslant n} \Range(X^i)) = 1\} \cap \{ |\Phi_{m'}| \leqslant 200 \} \Big) \bigg].
\end{split}
\end{equation*}

We now decompose the event under the above probability over all possible values of the number of excursions $G_i$ (performed by each of these $n$ trajectories) and on the departure points of these excursions ($X_{R_j}^i$ and $X_{D_j}^i$, for $j \leqslant G_i$). In the sums below, we tacitly assume that $x_j^i \in \partialint \Sigma$ and $y_j^i \in \partial C_{m(x_j^i)}^{L/4}$, for $i \leqslant n$ and $j \leqslant g_i$. The term above is thus smaller or equal to
\begin{equation}
\begin{split}
\leqslant 2^{2^\K} \sup_{\mbox{\fontsize{8}{8} \selectfont $\substack{M'\subset M; \\ |M'| = 2^{\K-1}}$}} \,\, \sum_{\mbox{\fontsize{8}{8} \selectfont $\substack{n \geqslant 0, g_1, \dots, g_n \geqslant 1\\ x_j^i, y_j^i \text{ such that } \forall m' \in M' \\ \#\{ (i,j); m(x_j^i) = m' \} \leqslant 200 }$}} \mathbb{P}[S=n] P_{\bar e_\Sigma}^{\otimes n} \bigg[ G_i = g_i, X_{R_j}^i = x_j^i, X_{D_j}^i = y_j^i,&\\
\qquad \text{ for $i \leqslant n$ and $j \leqslant g_i$}, \bigcap_{m'\in M'} \{ \chi_{m'}(\cup_{(i,j) \in \Phi_{m'}} E_j^i ) = 1\}\bigg]\\
\leqslant 2^{2^\K} \sup_{\mbox{\fontsize{8}{8} \selectfont $\substack{M'\subset M; \\ |M'| = 2^{\K-1}}$}} \,\, \sum_{\mbox{\fontsize{8}{8} \selectfont $\substack{n \geqslant 0, g_1, \dots, g_n \geqslant 1\\ x_j^i, y_j^i \text{ such that } \forall m' \in M' \\ \#\{ (i,j); m(x_j^i) = m' \} \leqslant 200 }$}} \mathbb{P}[S=n] E_{\bar e_\Sigma}^{\otimes n} \bigg[ G_i = g_i, X_{R_j}^i = x_j^i, X_{D_j}^i = y_j^i,&\\
\qquad \text{ for $i \leqslant n$ and $j \leqslant g_i$}, \prod_{m'\in M'} P_{\mbox{\fontsize{8}{8} \selectfont $\substack{(x_j^i)_{(i,j) \in \Phi_{m'}},(y_j^i)_{(i,j) \in \Phi_{m'}}}$}} \big[ \chi_{m'}(\cup_{(i,j) \in \Phi_{m'}} E_j^i ) = 1\big]\bigg].
\end{split}
\end{equation}
Choosing $G = 200$ and $\epsilon = (2^{10}(5 \cdot 80L)^{4d})^{-2}$ in Theorem~\ref{th:local}, we can bound the product above by $\epsilon^{2^{\K-1}}$ and remove it from the sum (which sum up to one), obtaining the desired bound
\begin{equation}
\label{eq:boundsecond}
\mathbb{P} \bigg[ \langle \mu_{\Sigma,u}, g_M \rangle \leqslant 100 \cdot 2^\K, \bigcap_{m\in M} \{ \chi_m(\mathcal{I}^u) = 1\}\bigg] \leqslant (2^{9}(5 \cdot 80L)^{4d})^{-2^\K},
\end{equation}
for an appropriate choice of $L_0 = L_0(L)$ which came from (\ref{th:reftolocal}).

With (\ref{eq:twoterms}), (\ref{eq:boundfirst}), (\ref{eq:etolambda}) and (\ref{eq:boundsecond}), we obtain
\begin{equation}
\label{eq:baruchoice}
\mathbb{P} \bigg[ \bigcap_{m\in M} \{ \chi_m(\mathcal{I}^u) = 1\}\bigg] \leqslant 2 \cdot (2^{9}(5 \cdot 80L)^{4d})^{-2^\K}
\end{equation}
for $u$ smaller or equal to some $\bar u = \bar u (L_0) > 0$. This, together with (\ref{th:coarse}) leads to (\ref{eq:reduct}), proving Theorem~\ref{th:main}. $\quad \square$

Note the order in which we chose the parameters for the proof: first $L = L(d)$ in (\ref{eq:Lchoice}), then $L_0 = L_0(L,d)$ in (\ref{eq:boundsecond}) and finally $\bar u = \bar u(L_0,L,d)$ in (\ref{eq:baruchoice}).

\vspace*{4mm}

\begin{remark}
\label{rem:G100}
$\,$

\textnormal{1) Let us comment on the constant $100$ appearing in the equation (\ref{eq:boundfirst}). In order to bound the first term of \eqref{eq:twoterms} we use an exponential Chebychev inequality, which provides a decay of $\mathbb{P}[\langle \mu_{\Sigma,u}, g_M \rangle > a \, 2^\K]$ of order $\sim \negmedspace (\text{const} \cdot L^{d-2})^{-a \cdot 2^\K}$, see \eqref{eq:boundfirst}. Notice that this decay should be fast enough in order to offset the growth of the combinatorial factor in \eqref{th:coarse}, $ii)$ (which is $\sim \negmedspace (\text{const}\cdot L^{4d})^{2^\K}$).}

\textnormal{The way we have to tune these two competing terms is by choosing $L$ large, but both of them depend on this parameter. In order to make this competition to work to our advantage, we need to choose a large number $a$ above, in such a way that $a(d-2) > 4d$. For our purposes $a = 100$ will do the job. This delicate balance is well illustrated in \eqref{eq:Lchoice} and \eqref{eq:boundfirst}, where we choose a large $L$.}


\textnormal{2) Note that the only part of the proof of Theorem~\ref{th:main} in which we use the hypothesis $d \geqslant 5$ is Theorem~\ref{th:local}. This will be further discussed in the Section~\ref{sec:local}.}
\end{remark}

\section{Coarse graining}
\label{sec:coarse}

In this section we study the hierarchical property of the function $\chi$, see Theorem~\ref{th:coarsefirst}, which was used to prove Theorem~\ref{th:main}.

\begin{figure}
\psfrag{A1}[cl][cl][3][0]{$A_1$}
\psfrag{A2}[cl][cl][3][0]{$A_2$}
\psfrag{U}[cl][cl][3][0]{$\U$}
\psfrag{C0}[cl][cl][3][0]{$C_{m(0)}$}
\psfrag{CT}[cl][cl][3][0]{$C_{m(N_\tau)}$}
\psfrag{CL}[cl][cl][3][0]{$C_{m(\bar l)}$}
\begin{center}
\includegraphics[angle=0, width=0.35\textwidth]{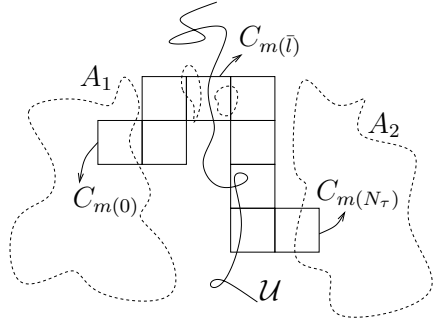}\\
\caption{The existence of $\bar l$, see Lemma~\ref{lem:indicespath}.}\label{fig:lbar}
\end{center}
\end{figure}

The main step to establish (\ref{th:coarse}) is the lemma below, which is interesting by itself, see the paragraph above (\ref{eq:reducetoK}). Loosely speaking it states that, if a path of boxes connects two large connected sets which are separated by $\U \subset \mathbb{Z}^d$, then at least one of this boxes (say indexed by $m'$) satisfies $\chi_{m'}(\U) = 1$, see Figure~\ref{fig:lbar}. More precisely,

\begin{lemma}
\label{lem:indicespath}
Let $A_1, A_2 \subset B \subset \mathbb{Z}^d$ be connected sets and $\U \subset \mathbb{Z}^d$ be such that $A_1$ is separated from $A_2$ by $\U$ in $B$. Consider some scale $\K \geqslant 0$ and a path $\tau$ in $\mathbb{Z}^d$, as well as the path of indices $m(l) = (\K, \tau(l)) \in I_\K$, $l = 0, \dots, N_\tau$, so that the following holds
\begin{equation}
\label{eq:indicespath}
\begin{aligned}
i) & \,\, C_{m(l)}^3 \subset B \text{ for all } l = 0,\dots, N_\tau,\\
ii) & \,\, C_{m(0)} \cap A_1 \neq \varnothing, C_{m(N_\tau)}\cap A_2 \neq \varnothing,\\
iii) & \,\, \diam(A_1), \diam(A_2) \geqslant L_\K/2.
\end{aligned}
\end{equation}
Then there is an $\bar l \in \{0,\dots, N_\tau\}$ such that $\chi_{m(\bar l)}(\U) =1$. See Figure~\ref{fig:lbar}.
\end{lemma}

\textit{Proof.} We consider the set
\begin{equation}
\label{eq:setA}
\begin{split}
A = \big\{ x\in B; & \text{ there is a path in $B$ from $\partial A_1$ to $x$}\\
& \text{ which is disjoint from $\U$} \big\} \cup A_1.
\end{split}
\end{equation}
We claim that
\begin{equation}
\label{eq:Aconnected}
A \text{ is connected}.
\end{equation}
Indeed, given $x,y \in A$, we join them to $\partial A_1$ by paths as in (\ref{eq:setA}) (and then to $A_1$ in case they were not already there) and use the connectedness of $A_1$.

Moreover
\begin{equation}
\label{eq:AfarfromA2}
A \cap \overline{A_2} = \varnothing.
\end{equation}
Indeed, suppose that there is some $x \in A\cap \overline{A_2}$. Since $d(A_1,A_2) > 1$ (see Definition~\ref{def:separate}) we have that $x \not\in A_1$. By (\ref{eq:setA}), there is a path in $B$ from $\partial A_1$ to $x$ which does not intersect $\U$. Possibly stopping this path at the first time it meets $\partial A_2$, we obtain a contradiction with the fact that $A_1$ is separated from $A_2$ by $\U$ in $B$.

We then claim that
\begin{equation}
\label{eq:boundAinU}
\partial A \cap B \subset \U.
\end{equation}
Suppose by contradiction the existence of some $x \in (\partial A \cap B)\setminus \U$, say with $x$ neighbor of $y \in A$. If $y \in A_1$, we have $x \in (\partial A_1 \cap B) \setminus \U$ which implies by (\ref{eq:setA}) that $x\in A$, a contradiction. The other possibility is that $y$ is connected to $\partial A_1$ by some path in $B$ which is disjoint from $\U$. Adding one step to this path we can assume that it meets $x$, obtaining again the contradiction $x \in A$, see (\ref{eq:setA}). This proves (\ref{eq:boundAinU}).

In order to choose $\bar l$, we define
\begin{equation}
\label{eq:setI}
I = \{ l \in \{ 0,\dots , N_\tau\}; C_{m(l)} \cap A = \varnothing \}
\end{equation}
and consider the two following possibilities

\vspace*{4mm}

\textit{Case 1: $I$ is empty.}

In other words, the set $A$ meets all the boxes with indexes in the path $\tau$. In this case we take $\bar l = N_\tau$.

From (\ref{eq:setI}), there is some point $x \in C_{m(\bar l)} \cap A$ and by (\ref{eq:indicespath}) $ii)$, there is some $y \in C_{m(\bar l)}\cap A_2$. By (\ref{eq:AfarfromA2}) and (\ref{eq:boundAinU}), $A$ and $A_2$ are separated by $\U$ in $B$. Since (\ref{eq:indicespath}) $i)$ implies that all paths in $C_{m(\bar l)}^3$ are in $B$, all we need in order to show that $\chi_{m(\bar l)}(\U) = 1$ is to extract connected components $A_1'\subset A \cap C_{m(\bar l)}^2, A_2' \subset A_2 \cap C_{m(\bar l)}^2$ with $\diam(A_1'), \diam(A_2') \geqslant \tfrac{L_\K}{2}$.

We know that
\begin{equation}
\label{eq:diamA}
\diam(A) \geqslant \diam (A_1) \geqslant L_\K /2.
\end{equation} 
If $ A \subset C_{m(\bar l)}^2$, since $A$ is itself connected, we take $A_1' = A$. Otherwise, again using the connectedness of $A$, we take $A_1'$ to be the range of a path from $x$ to $\partialint C_{m(\bar l)}^2$.

With a similar argument we obtain $A_2'$, since it is connected, intersects $C_{m(\bar l)}$ and its diameter is at least $L_\K/2$.

\vspace*{4mm}

\textit{Case 2: $I$ is non-empty.}

In this case, we have by (\ref{eq:indicespath}) $ii)$, that $0 \not\in I$, so that $\min \{I\} > 0$ and we define $\bar l$ to be $(\min \{I\} -1)$. We know that $C_{m(\bar l)} \cap A \neq \varnothing$ and as in the previous case, we are able to find some connected set $A_1' \subset A \cap C_{m(\bar l)}^2$ such that $\diam(A_1') \geqslant L_\K/2$. Define $A_2'$ to be $C_{m(\bar l +1)} \setminus \partialint C_{m(\bar l + 1)}$, which is connected and contained in $C_{m(\bar l)}^2$. By the definition of $\bar l$, we have that $C_{m(\bar l +1)}\cap A = \varnothing$, hence, $d(A_1', A_2') \geqslant d(A,A_2') >1$. From (\ref{eq:boundAinU}) we conclude that any path in $B$ from $A_1'$ to $A_2'$ meets $\partial A$ and hence $\U$. And since $\diam(A_2') \geqslant L_\K/2$, we find that $\chi_{m(\bar l)}(\U) = 1$.

This finishes the proof of Lemma~\ref{lem:indicespath}. $\quad \square$

\vspace*{4mm}

Now we prove the following
\begin{theorem}
\label{th:coarseinduct}
Let $m \in I_{\K +1}$ and $\U \subset \mathbb{Z}^d$. If $\chi_m(\U) = 1$, we can find $m_1 = (\K,i_1)$ and $m_2 = (\K,i_2)$ such that
\begin{equation}
\label{eq:coarseinduct}
\begin{aligned}
i) & \,\, C_{m_1}^3, C_{m_2}^3 \subset C_m^3,\\
ii) & \,\, |i_1 -i_2| \geqslant 2L,\\
iii)& \,\, \chi_{m_1}(\U) = \chi_{m_2}(\U) = 1.\\
\end{aligned}
\end{equation}
\end{theorem}

\textit{Proof.} Since $\chi_m(\U) = 1$, we can find connected sets $A_1, A_2 \subset C_m^2$ such that $\diam(A_1)$, $\diam(A_2)$ $\geqslant L_{\K+1}/2 > L_\K/2$ and $A_1$, $A_2$ are separated by $\U$ in $C^3_m$.

With the Lemma~\ref{lem:indicespath}, we conclude that Theorem~\ref{th:coarseinduct} follows once we show that
\begin{equation}
\label{eq:allweneed}
\begin{array}{c}
\text{there are two paths $\tau$ and $\tau'$ such that $d_\infty(\Range(\tau),\Range(\tau')) \geqslant 2L$ and the}\\
\text{corresponding paths of incices at level $\K$ satisfy the conditions (\ref{eq:indicespath}) $i)$ and $ii)$.}
\end{array}
\end{equation}

Since $\diam(A_i) \geqslant L_{\K + 1} / 2$ ($i = 1,2$), we can find
\begin{equation}
\label{eq:i1i2andprime}
\begin{array}{c}
i_1,i_2,i_1',i_2' \in \mathbb{Z}^d \text{ such that }
C_{(\K,i_j)}\cap A_j, C_{(\K,i_j')} \cap A_j \neq \varnothing\\
\text{and } |i_j - i_j'| \geqslant 20L, \text{ for } j = 1, 2.
\end{array}
\end{equation}
After relabeling, we can suppose that
\begin{equation}
\label{eq:noprimeismin}
\text{the pair } i_1,i_2 \text{ minimizes the $\infty$-distance between the sets } \{i_1,i_1'\} \text{ and } \{i_2,i_2'\},
\end{equation}
and we set $\Lambda = \lVert i_1 - i_2 \rVert_\infty$. We claim that
\begin{equation}
\label{eq:distprimeandno}
\lVert i_1' - i_2 \rVert_\infty, \lVert i_2' - i_1 \rVert_\infty \geqslant 10L.
\end{equation} 
Indeed, with (\ref{eq:noprimeismin}) and the triangle inequality,
\begin{equation}
\begin{split}
\lVert i_1' - i_2 \rVert_\infty & \geqslant \max \{ \Lambda, \lVert i_1 - i_1' \rVert_\infty - \Lambda\}\\
& \geqslant \max \{\Lambda, 20L- \Lambda \} \geqslant 10L.
\end{split}
\end{equation}

As a result of (\ref{eq:i1i2andprime}) and (\ref{eq:distprimeandno}), we find that
\begin{equation}
\label{eq:distprimetoi1}
d_\infty(\{i_1',i_2'\},\{i_1, i_2\}) \geqslant 10L.
\end{equation} 

To find the desired paths $\tau$ and $\tau'$, we consider two cases:

\begin{figure}
\psfrag{i1}{$i_1$}
\psfrag{j1}{$i_1$}
\psfrag{i1p}{$i_1'$}
\psfrag{j1p}{$i_1'$}
\psfrag{i2}{$i_2$}
\psfrag{j2}{$i_2$}
\psfrag{i2p}{$i_2'$}
\psfrag{j2p}{$i_2'$}
\psfrag{b1}{$B(i_1,8L)$}
\psfrag{b2}{$B(i_1,10L)$}
\psfrag{b3}{$B(i_1,6L)$}
\psfrag{b4}{$B(i_2,6L)$}
\psfrag{t}{$\tau$}
\psfrag{tp}{$\tau'$}
\psfrag{t2}{$\tau$}
\psfrag{tp2}{$\tau'$}
\psfrag{F}{$F$}
\psfrag{D}{$D$}
\psfrag{b}{$\partial^* B(i_1,10L)$}
\psfrag{bp}{$\partial^* D$}
\begin{center}
\includegraphics[angle=0, width=1.0\textwidth]{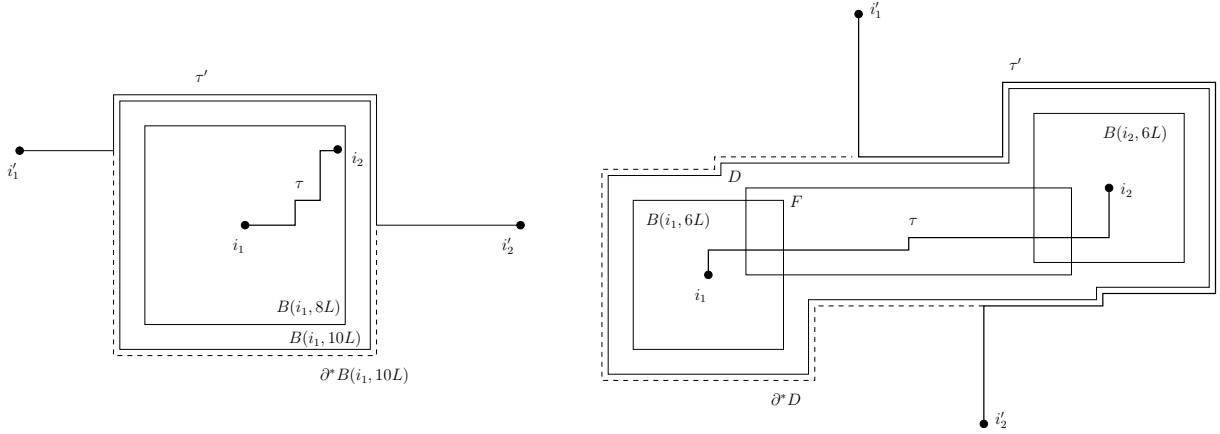}\\
\caption{The choice of $\tau$ and $\tau'$ in \textit{Case 1.} (left) and \textit{Case 2.} (right)}\label{fig:paths}
\end{center}
\end{figure}

\vspace*{4mm}
\textit{Case 1:} $\Lambda < 8L$.

In this case we connect $i_1$ to $i_2$ by any path $\tau$ contained in $B(i_1,8L)$, see Figure~\ref{fig:paths}.

In order to find $\tau '$, we first connect $i_1'$ to $i_2'$ by any $\bar \tau$ in $\text{sbox}(\{i_1', i_2'\})$. Since $C_{(\K,i_1)},$ $C_{(\K,i_2)} \subset C_m^2$, we have
\begin{equation}
C_{(\K,\bar \tau(n))} \subset C_m^2, \text{ for all } n=0,\dots, N_{\bar \tau}.
\end{equation}

By (\ref{eq:distprimetoi1}) and the connectedness of $\partial^* B(i_1,10L)$, see (\ref{eq:stfillconn}), we can find a modification $\tau'$ of the path $\bar \tau$ which joins the points $i_1'$ and $i_2'$, avoiding $B(i_1,10L)$ and satisfying $\Range(\tau') \subset B(\Range(\bar \tau), 20L + 2)$. So that $d_\infty(C_{(\K,\tau'(n))},C_m^2) \leqslant 30L \cdot L_{\K} \leqslant L_{\K +1}/2$ and consequently $C^3_{(\K,\tau'(n))} \subset C_m^3$ for $n=0, \dots, N_{\tau'}$. Finally, $\Range(\tau) \subset B(i_1,8L)$ and $\Range(\tau') \subset B(i_1,10L)^c$ imply that $d_\infty(\Range(\tau),\Range(\tau')) \geqslant 2L$.

\vspace*{4mm}

\textit{Case 2:} $\Lambda \geqslant 8L$.

First we claim that
\begin{equation}
\label{eq:setF}
F = B(i_1,\Lambda -3L) \cap B(i_2,\Lambda -3L) \cap \text{sbox}(\{i_1,i_2\}) \text{ is a non-empty box},
\end{equation}
see Figure~\ref{fig:paths}. Indeed, if we take $y = \text{floor}(\frac{i_1+i_2}{2})$ (recall the definition below (\ref{eq:floorconvex})), we have
\begin{enumerate}
\item $\lVert y - i_j \rVert_\infty \leqslant \tfrac{\Lambda}{2} + 1 \leqslant \Lambda -3L$ for $j=1,2$, since $\Lambda \geqslant 8L$,
\item $y \in \text{sbox}(\{i_1,i_2\})$, applying to each coordinate the equation (\ref{eq:floorconvex}).
\end{enumerate}
Since, a non-empty intersection of boxes is a box, (\ref{eq:setF}) follows.

We also claim that
\begin{equation}
\label{eq:Fcapballs}
\text{$F$ intersects $B(i_1,6L)$ and $B(i_2,6L)$}.
\end{equation}
For this, take $y = \text{floor}\big(\tfrac{(3L)i_2+(\Lambda -3L)i_1}{\Lambda}\big)$, using once more (\ref{eq:floorconvex}) we conclude that $y \in \text{sbox}(\{i_1,i_2\})$. We then bound $\lVert y - i_1 \rVert_\infty$ as follows:
\begin{equation}
\lVert y - i_1 \rVert_\infty \leqslant 1+ \Big\Vert \tfrac{(3L)i_2+(\Lambda -3L)i_1}{\Lambda} - i_1 \Big\Vert_\infty \leqslant 1 + \tfrac{3L}{\Lambda} \lVert i_2 -i_1 \rVert_\infty \leqslant 4L.
\end{equation}
We thus find that $y \in F \cap B(i_1,6L)$. In a similar manner we see that $y \in F \cap B(i_2,6L)$, finishing the proof of the claim (\ref{eq:Fcapballs}).

We now connect $i_1$ to $i_2$ by any path $\tau$ in $F \cup B(i_1,6L) \cup B(i_2,6L) \subset B(F,12L)$. Since $d_\infty(C_m^2, (C_m^3)^c) = L_{\K +1} > 20 L \cdot L_\K$, we have $C_{(\K,\tau(n))}^3 \subset C_m^3$ for all $n = 1, \dots, N_\tau$.

Finally, define $D = B(F,2L) \cup B(i_1,8L) \cup B(i_2,8L)$. Equation (\ref{eq:distprimetoi1}) implies that $\{i_1', i_2'\} \cap (B(i_1,8L)\cup B(i_2,8L)) = \varnothing$. By (\ref{eq:noprimeismin}), we have $\lVert i_1' - i_2 \rVert_\infty , \lVert i_2' - i_1 \rVert_\infty \geqslant \lVert i_2 - i_1 \rVert_\infty = \Lambda$. Hence, $\{i_1',i_2'\} \cap B(i_1,\Lambda -L) \cap B(i_2,\Lambda -L) = \varnothing$ and neither $i_1'$ nor $i_2'$ belongs to $D$.

Take $y$ in $D^c$. Since only three of the $2d$ half-lines (parallel to the canonical basis) connecting $y$ to infinity can meet $D$ ($D$ is the union of three boxes), we see with Lemma~\ref{eq:charfill} that $y \not \in \f(D)$. So that $\f(D) = D$ and (\ref{eq:stfillconn}) implies that $\partial^* D$ is connected. As in the \textit{Case 1}, we choose some path $\bar \tau$ in $\text{sbox}(\{i_1',i_2'\})$, connecting $i_1'$ to $i_2'$ and modify it to get some $\tau'$ joining $i_1'$ to $i_2'$ which is disjoint from $D$ and contained in $\text{sbox}(\{i_1',i_2'\}) \cup \overline{D} \subset C_m^3$.

Together with (\ref{eq:allweneed}), this concludes the proof of Theorem~\ref{th:coarseinduct}. $\quad \square$

\vspace*{4mm}

In Section~\ref{sec:proof}, we stated Theorem~\ref{th:coarsefirst}, which we now deduce from the Theorem~\ref{th:coarseinduct}. Note that the statement of Theorem~\ref{th:coarse2} below is slightly stronger than Theorem~\ref{th:coarsefirst} (as needed in the induction procedure used in the proof) and thus the following implies Theorem~\ref{th:coarsefirst}.

\begin{theorem}
\label{th:coarse2}
Given $m' \in I_\K$, $\K \geqslant 0$, there is a family of skeletons $\mathcal{M}_{m'}$, see Definition~\ref{def:skeleton}, such that
\begin{equation}
\label{eq:coarse2}
\begin{aligned}
i) &\,\, \text{for all $M \in \mathcal{M}_{m'}$, $\#M = 2^{\K}$ and for $m \in M$, $C^3_m \subset C^3_{m'}$,}\\
ii) &\,\, \text{$\#\mathcal{M}_{m'} \leqslant ((5 \cdot 80L)^{4d})^{2^\K-1}$,}\\
iii) &\,\, \chi_{m'}(\U) \leqslant \sum_{M \in \mathcal{M}_{m'}} \prod_{m \in M} \chi_m(\U), \text{ for all $\U \subset \mathbb{Z}^d$, ``if $\U$ promotes a separation}\\
& \,\, \text{in $C^3_{m'}$, it also separates components in all boxes of a skeleton $M \in \mathcal{M}_{m'}$''}
\end{aligned}
\end{equation}
\end{theorem}

\textit{Proof.} We proceed by induction on $\K$. If $\K =0$, (\ref{eq:coarse2}) holds for $\mathcal{M}_{m'} = \{\{m'\}\}$.
Assume now that (\ref{eq:coarse2}) holds for some $\K$ and consider for $m'\in I_{\K+1}$,
\begin{equation}
\label{eq:Mmprimeequals} 
\mathcal{M}_{m'} = \bigcup_{\substack{m_1 = (\K,i_1), m_2 = (\K,i_2); \\ C^3_{m_1},C^3_{m_2} \subset C^3_{m'}, |i_1 -i_2| \geqslant 2L }} \big\{M_1 \cup M_2; M_1 \in \mathcal{M}_{m_1}, M_2\in \mathcal{M}_{m_2}\big\}.
\end{equation}

Note that (\ref{eq:coarse2}) $iii)$ directly follows from Theorem~\ref{th:coarseinduct} and the induction hypothesis.

Next we will show that $\mathcal{M}_{m'}$ is in fact a family of skeletons, see Definition~\ref{def:skeleton}. To this end, consider $m_0 \in M_1 \cup M_2$, for $M_1, M_2$ as in (\ref{eq:Mmprimeequals}). We suppose without loss of generality that $m_0 \in M_1$. By the induction hypothesis, $M_1$ is a skeleton, so that
\begin{gather}
\label{eq:M1skeleton1}
\# \{m\in M_1; L\,L_h \leqslant d_\infty(C_{m_0},C_m) \leqslant L\,L_{h+1}\} \leqslant 2^{h+1},\\
\label{eq:M1skeleton2}
\{m\in M_1; d_\infty(C_{m_0},C_m) \leqslant L\,L_0\} = \{m_0\}.
\end{gather}

Since all boxes $C_m$ for $m \in M_1$ are contained in $C^3_{m_1}$ ($\subset C^3_{m'}$) and $\diam(C^3_{m_1}) = 5 L_\K < L\, L_\K$, the sets in (\ref{eq:M1skeleton1}) are empty for $h \geqslant \K$.

From the inequality $|i_1 -i_2| \geqslant 2L$, we deduce that
\begin{equation}
\label{eq:distM1M2}
\text{for every } m\in M_2, d_\infty(C_m,C_{m_0}) \geqslant d_\infty(C^3_{m_1},C^3_{m_2}) \geqslant (2L-5) L_\K > L\,L_\K.
\end{equation}
Hence,
\begin{multline}
\#\{ m \in M_1 \cup M_2; L\,L_h < d_\infty(C_{m_0},C_m)\leqslant L\,L_{h+1}\}\\
= \begin{cases}
\{m \in M_1;L\,L_h < d_\infty(C_{m_0},C_m)\leqslant L\,L_{h+1} \} \leqslant 2^{h+1}& \text{if $h < \K$},\\
\{m \in M_2;L\,L_h < d_\infty(C_{m_0},C_m)\leqslant L\,L_{h+1} \} \leqslant \# M_2 = 2^\K \leqslant 2^{h+1} & \text{if $h \geqslant \K$}.
\end{cases}
\end{multline}
and
\begin{equation}
\{ m \in M_1 \cup M_2; d_\infty(C_{m_0},C_m) \leqslant L\,L_{0}\} = \{ m \in M_1; d_\infty(C_{m_0},C_m) \leqslant L\,L_{0}\} = \{m_0\}.
\end{equation}

This shows that $M_1 \cup M_2$ is a skeleton.

By (\ref{eq:distM1M2}), $M_1$ and $M_2$ are disjoint, so that $\# (M_1 \cup M_2) = 2^{\K +1}$ and (\ref{eq:coarse2}) $i)$ holds. Moreover,
\begin{equation}
\begin{split}
\# \mathcal{M}_{m'} & \leqslant \# \{(m_1,m_2); m_1 \in I_\K, m_2 \in I_\K \text{ and } C^3_{m_1}, C^3_{m_2} \subset C^3_{m'}\} \cdot \#\mathcal{M}_{m_1} \cdot \#\mathcal{M}_{m_2}\\
& \leqslant (5 \cdot 80L)^{2d}\Big(\big((5 \cdot 80L)^{4d}\big)^{2^\K -1}\Big)^2\\
& = (5 \cdot 80L)^{2d}\big((5 \cdot 80L)^{4d}\big)^{2^{\K + 1}-2} \leqslant \big( (5 \cdot 80L)^{4d} \big)^{2^{\K+1} - 1},
\end{split}
\end{equation}
and (\ref{eq:coarse2}) $ii)$ is verified. This concludes the proof by induction of Theorem~\ref{th:coarse2}. $\quad \square$

\section{Walking around sausages}
\label{sec:local}

The aim of this section is to establish Theorem~\ref{th:local} (already stated in Section~\ref{sec:proof} and used in the proof of our main result, Theorem~\ref{th:main}, above (\ref{eq:boundsecond}).

Roughly speaking, Theorem~\ref{th:local} states that with overwhelming probability, a fixed number of random walk trajectories does not separate components of macroscopic diameter in a large enough box. Moreover, the statement of Theorem~\ref{th:local} holds uniformly over the points at which we condition these random walk trajectories to enter and exit large neighborhoods of this box.

Actually, the current section is the only part of the proof of Theorem~\ref{th:main} where the hypothesis ($d \geqslant 5$) is used, so that in order to extend Theorem~\ref{th:main} to lower dimensions, it would be enough to prove a version of Theorem~\ref{th:local} for this case, see also Remark~\ref{rem:d3}.

\vspace*{4mm}

We now give a rough overview of the proof of Theorem~\ref{th:local} which relies on Lemmas~\ref{lem:conditioned}, \ref{lem:farapart} and \ref{lem:Gfarapart}.

The first step of the proof consists in proving an analogue of Theorem~\ref{th:local} for one single trajectory. More precisely, in Lemma~\ref{lem:conditioned} we prove that the probability that one random walk trajectory separates components of macroscopic diameter inside a box goes to zero as the size of this box increases. Moreover, this limit is uniform over the points $x$ and $y$ at which we condition this random walk to enter and exit large neighborhoods of the box.

To prove this lemma, we regard a `chunk' $\U$ of the random walk trajectory as a set of `sausages' connected by cut-points, see (\ref{eq:cuttime}). An important concept here is the notion of $h$-avoidable sets, where $h \geqslant 1$ is an integer, see Definition~\ref{def:avoidable}. Loosely speaking, a set $A$ is said to be $h$-avoidable if any path traversing $A$ can be modified (within a distance of at most $h$) in order to go around $A$ through its boundary. In Remark~\ref{rem:B01} we exemplify this definition showing that $B(0,1)$ is $4$-avoidable, while $\partialint B(0,1)$ is not.

The heart of the proof of Lemma~\ref{lem:conditioned} is Lemma~\ref{lem:avoid}, which roughly states the following: for a piece of trajectory $\U$, if the diameters of its `sausages' are bounded by $h$, then $\f(\U)$ is $(3h)$-avoidable. The strategy to prove this lemma can be informally described as ``to travel through the skins of the sausages''. This proof clarifies and solves the geometric restrictions mentioned in Remark~\ref{rem:Bernoulli} 3).

As a direct consequence of Lemma~\ref{lem:avoid}, we conclude that for any pair of sets $A_1$ and $A_2$, which are large when compared to the diameter of the `sausages' in $\U$, we can connect $\partial A_1$ to $\partial A_2$ avoiding $\U$. In Corollary~\ref{lem:chifillX}, we conclude that if the diameters of these `sausages' are bounded by $L_0/4$, then $\U$ does not separate components in a box of diameter $L_0-1$ (in the sense of (\ref{eq:chi})).

For $d \geqslant 5$, we are able to bound the diameter of the `sausages' occurring on a typical random walk trajectory before it exits the neighborhood of a given box, see Lemma~\ref{lem:estimates}. The proof of Lemma~\ref{lem:estimates} relies on known results on intersections of random walks, see for instance \cite{lawler}.

The above mentioned results conclude the proof of Lemma~\ref{lem:conditioned}. The uniformity of this lemma over the points $x$ and $y$ (where we condition the random walk to enter and exit large neighborhoods of the given box) follows from the Harnack inequality, see (\ref{eq:Harnack}).

\vspace*{4mm}

The second step in the proof of Theorem~\ref{th:local} is to extend Lemma~\ref{lem:conditioned} from one trajectory to a fixed number, say $G$, of independent random walk trajectories. First, consider $G$ connected subsets of $\mathbb{Z}^d$ such that: none of them separates components in a certain box, they are mutually far apart and they are $(L_0/2G)$-avoidable. In Lemma~\ref{lem:farapart} we show that the union of these $G$ sets also does not separate components in the box.

Thus, all that remains to prove is that $G$ independent random walks are, with high probability, mutually far apart. This is the content of Lemma~\ref{lem:Gfarapart}, which again uses arguments on intersections of random walks for $d \geqslant 5$.

Finally, we bring together Lemmas~\ref{lem:conditioned}, \ref{lem:farapart} and \ref{lem:Gfarapart} to obtain Theorem~\ref{th:local}.

\vspace*{4mm}

Let us introduce the concept of cut-times for a doubly infinite trajectory. Let $w \in W$ and recall that $\{X_i\}_{i \in \mathbb{Z}}$ denote its canonical coordinates.
\begin{equation}
\label{eq:cuttime}
\text{We say that $k \in \mathbb{Z}$ is a \textit{cut-time} of $w$ if $d_\infty(\{X_l\}_{l<k}, \{X_l\}_{l>k}) > 1$}.
\end{equation}
In this case $X_k$ is called a \textit{cut-point}. Note that our definition differs from the usual definition of a cut-time, which does not require a strict inequality as above, see for instance \cite{erdos}.

What we informally described as a 'sausage' will be determined by the range of $w$ between two chosen cut-times. Note, however, that the definition in (\ref{eq:cuttime}) does not exclude the possibility that two cut-times are adjacent (e.g. for the trajectory $X_j = j{\pmb e}_1$, every integer is a cut-time). So, given a finite sequence of cut-times $n_0 < \dots < n_J$ which are not adjacent, i.e.
\begin{equation}
\label{eq:farcuts}
n_{j+1} - n_j \geqslant 2, \text{ for all } j = 0, \dots, J-1,
\end{equation}
we define the `pieces of trajectory':
\begin{equation}
\label{eq:Uj}
\U_j = X_{(n_j,n_{j+1})}, \text{ for } j = 0, \dots, J-1.
\end{equation}
We stress here that the definition of sausages depend on the choice of the (non-adjacent) cut-times, which in general will not be the whole set of cut-times in a given interval. There is still another technical reason for the introduction of the condition (\ref{eq:farcuts}), see the proof of Lemma~\ref{lem:avoid}, above (\ref{eq:Utilde}). The sets $\f(\U_j)$ are what we informally referred to as 'sausages', recall Definition~\ref{def:fill}.

\vspace*{4mm}

The next lemma states some useful properties of $\f(A)$ which we exclusively need in this section. Its proof, although short, digresses from our main purpose here and can be found in the Appendix.

\begin{lemma}
\label{prop:fill}
If $A \subset \mathbb{Z}^d$ is finite, then
\begin{gather}
\label{eq:boundfill}
\partialint \f(A) \subset A \text{ and}\\
\label{eq:diamfill}
\diam(\f(A)) = \diam(A).
\end{gather}
Now let $A,B \subset \mathbb{Z}^d$ be connected and finite. If $d(A,B) > 1$ (respectively $d_\infty(A,B) > 1$) then exactly one of the following possibilities holds:
\begin{equation}
\label{eq:threecases}
\begin{aligned}
i)&\,\, A \subset \f(B) \text{ and } d(\f(A),B) > 1, \text{``$A$ is interior to $B$''},\\
ii)&\,\, B \subset \f(A) \text{ and } d(\f(B),A) > 1, \text{``$B$ is interior to $A$''}, \\
iii)&\,\,d(\f(A),\f(B)) > 1 \text{ (respectively $d_\infty(\f(A),\f(B)) > 1$)}, \\
&\,\, \text{``$A$ and $B$ are exterior to each other''}.
\end{aligned}
\end{equation}
\end{lemma}

\vspace*{4mm}

As we described at the beginning of this section, our main argument to show that certain sets do not separate components relies on the definition of `avoidable set'. This is made precise in the following

\begin{definition}
\label{def:avoidable}
For $h \geqslant 1$ and $A,C \subset \mathbb{Z}^d$ with $A$ finite, we say that $A$ is $h$-avoidable in $C$ if for every path $\tau$ in $C$ with endpoints $\tau(0), \tau(N_{\tau})$ not in $A$ we can find a modification $\tau'$ of $\tau$ such that:
\begin{equation}
\begin{aligned}
i)& \,\, \tau' \text{ does not meet $A$},\\
ii)&\,\, \tau'(0) = \tau(0) \text{ and } \tau'(N_{\tau'}) = \tau(N_\tau), \text{``$\tau$ and $\tau'$ have the same endpoints''},\\
iii) &\,\, \Range(\tau') \subset B(\Range(\tau),h) \cap (\Range(\tau) \cup \partial^* A), \text{``$\tau'$ remains $h$-close} \\
&\,\, \text{to $\Range(\tau)$ and outside $\Range(\tau)$ it stays in $\partial ^* A$''}.
\end{aligned}
\end{equation}
When $C = \mathbb{Z}^d$, we simply write that $A$ is $h$-avoidable. And if the value of $h$ is not relevant, we omit it in the notation.
\end{definition}
Note that we do not require the path $\tau'$ or set $A$ to be contained in $C$. However,
\begin{equation}
\label{eq:avoidremark}
\text{the property ``$A$ is $h$-avoidable in $C$'' only depends on the set $A \cap B(C,h)$.}
\end{equation}

\begin{figure}
\psfrag{a}{$A$}
\psfrag{ba}{$\partial ^* A$}
\psfrag{t}{$\tau$}
\psfrag{a2}{$\partialint A$}
\begin{center}
\includegraphics[angle=0, width=0.75\textwidth]{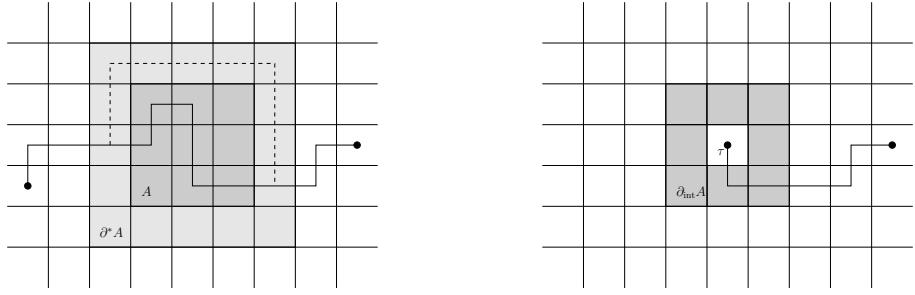}\\
\caption{The set $A = B(0,1)$ is $4$-avoidable, however, $\partialint A \subset A$ is not.}\label{fig:B01}
\end{center}
\end{figure}

\begin{remark}
\label{rem:B01}
\textnormal{The property of being $h$-avoidable is not monotonic. Consider for instance the set $A = B(0,1)$ and some path with endpoints in $A^c$, see Figure~\ref{fig:B01}. Every excursion this path performs inside the set $A \cup \partial ^* A$ can be replaced by an excursion entirely contained in $\partial^* A$ (according to (\ref{eq:stfillconn}), $\partial ^* A = \partial^* \f(A)$ is connected). Since $\diam(\partial ^* A) = 4$, we conclude that $A$ is $4$-avoidable. Although $A$ is $4$-avoidable and the set $\partialint A$ is contained in $A$, we check that $\partialint A$ fails to be avoidable. Indeed, no path $\tau$ connecting the origin to some point in $A^c$ can be modified to another path which is disjoint from $\partialint A$ but have the same endpoints as $\tau$, see Figure~\ref{fig:B01}.}

\textnormal{Note also that,}
\begin{equation}
\label{eq:avoidmonotone}
\textnormal{if $A$ is $h$-avoidable in $C$ and $C' \subset C$, $A$ is also $h$-avoidable in $C'$.} \quad \square
\end{equation} 
\end{remark}

\vspace{4mm}

In the lemma below, we show that $\f(X_{(n_0,n_J)})$ is $(3h)$-avoidable, where $h$ is a bound on the diameter of the sets $\U_j$. Loosely speaking, we first show that the union of the `sausages' ($\f(\U_j)$) and the cut-points ($\{x_{n_j}\}$) is $(3h)$-avoidable, by ``traveling through the skins of the sausages''. Then we show that $\f(X_{(n_0,n_J)})$ is in fact precisely this union.

\begin{lemma}
\label{lem:avoid}
($d \geqslant 3$) If $n_0,\dots,n_J$ are cut-times chosen as in (\ref{eq:farcuts}) and for some $h \geqslant 2$,
\begin{equation}
\label{eq:diamcondit}
\max_{j=0,\dots, J-1} \{\diam(\U_j)\} \leqslant h,
\end{equation}
then $\f(X_{(n_0,n_J)})$ is $(3h)$-avoidable.
\end{lemma}

\vspace*{4mm}

\textit{Proof.} First we describe the neighborhood of the cut-points. Note that
\begin{equation}
\label{eq:cutaxis}
\underset{e}{\underbrace{X_{n_j + 1} - X_{n_j}}} = \underset{f}{\underbrace{X_{n_j} - X_{n_j - 1}}}, \text{ for $j = 1, \dots, J-1$}.
\end{equation}
Indeed, $\Vert e + f \Vert_\infty = \lVert X_{n_j + 1} - X_{n_j - 1} \rVert_\infty \overset{(\ref{eq:cuttime})}{>} 1$, and since $\lVert e \rVert = \lVert f \rVert = 1$, $e$ and $f$ must be equal.

\begin{figure}
\psfrag{a}{$X_{n_j-1}$}
\psfrag{b}{$X_{n_j}$}
\psfrag{c}{$X_{n_j+1}$}
\psfrag{d}{$e=f$}
\psfrag{e}{$Ring_j$}
\begin{center}
\includegraphics[angle=0, width=0.2\textwidth]{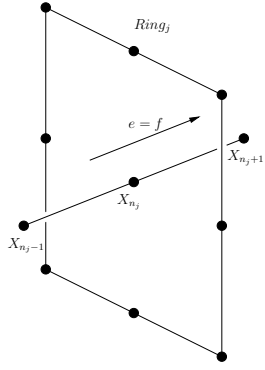}\\
\caption{The $Ring_j$ and the perpendicular vector $e$.}\label{fig:ring}
\end{center}
\end{figure}

We now define $Ring_j$, as the set of $\ast$-neighbors of $X_{n_j}$ which lay in the $(d-1)$-dimensional plane by $X_{n_j}$ perpendicular to $e$, i.e.
\begin{equation}
\label{eq:ringj}
Ring_j = \{ y \overset{*}{\leftrightarrow} X_{n_j}; (y-X_{n_j}) \perp e \},
\end{equation}
see Figure~\ref{fig:ring}.

We claim that for any $j = 1, \dots , J-1$,
\begin{equation}
\label{eq:ringdisj}
\text{$Ring_j$ is disjoint from the whole trajectory $X_{(-\infty,\infty)}$}.
\end{equation}
In fact, the definition of $Ring_j$ in (\ref{eq:ringj}) implies that
\begin{equation}
\label{eq:ringdist}
\text{every point $y$ in $Ring_j$ satisfies } \lVert y - X_{n_j - 1} \rVert_\infty = \lVert y - X_{n_j + 1} \rVert_\infty = 1,
\end{equation}
so that by (\ref{eq:cuttime}), $y$ is disjoint from both $\{ X_l \}_{l < n_j}$ and $\{ X_l \}_{l > n_j}$. Since $y \neq X_{n_j}$, this proves (\ref{eq:ringdisj}).

We also claim that the `sausages' are exterior to each other, i.e. for $k \geqslant 1$,
\begin{equation}
\label{eq:distsauss}
d_\infty(\f(\U_j),\f(\U_{j+k})) > 1, \text{ for } 0 \leqslant j < j+k \leqslant J-1.
\end{equation}
Indeed, by the definition of cut-times, the trajectory $(X_{n_{j+k} + 1 + i})_{i \geqslant 0}$ doesn't meet the set $\U_j$ (analogously, $(X_{n_{j+1} - 1 -i})_{i \geqslant 0}$ doesn't meet $\U_{j+k}$). So, using the characterization of $\f(\cdot)$ given in Lemma~\ref{eq:charfill}, we conclude that $\U_j \not \subset \f(\U_{j+k})$ and $\U_{j+k} \not \subset \f(\U_{j})$. By the the definition of the cut-time $n_{j+1}$, we have $d_\infty(\U_j,\U_{j+k}) > 1$, and using (\ref{eq:threecases}) of Lemma~\ref{prop:fill}, we obtain (\ref{eq:distsauss}).

As a consequence of (\ref{eq:ringdist}) and (\ref{eq:distsauss}), we have that
\begin{equation}
\begin{array}{c}
\label{eq:ringskin}
Ring_j \subset \partial^*\f(\U_{j-1}) \cap \partial^* \f(\U_{j}), \text{ for every } j =1,\dots, J-1, \\
\text{``the ring is contained in the skin of both of its adjacent sausages''}.
\end{array}
\end{equation}
And from (\ref{eq:stfillconn}) we conclude that
\begin{equation}
\label{eq:skinconn}
\text{$\partial^* \f(\U_j)$ is connected ($j = 0, \dots, J-1$).}
\end{equation}

We now show that the union of the `sausages' with the cut-points
\begin{equation}
\label{eq:Utilde}
\widetilde{\U} = \left[ \bigcup_{j=0}^{J-1} \f(\U_j) \right] \cup \{X_{n_j}\}_{j=1,\dots,J-1}
\end{equation}
is $(3h)$-avoidable and it will only remain to prove that $\f(X_{(n_0,n_{J})}) = \widetilde{\U}$.

As in the Definition~\ref{def:avoidable} of a $(3h)$-avoidable set, take a path $\tau$ such that
\begin{equation}
\label{eq:tauextr}
\tau(0), \tau(N_\tau) \notin \widetilde{\U}.
\end{equation}

Loosely speaking, we first modify $\tau$ into a path $\bar \tau$ which ``surrounds each sausage through its skin'' and finally we will modify $\bar \tau$ to a path $\tau'$ surrounding the cut-points using the rings.

Using (\ref{eq:distsauss}), we conclude that the visits performed by $\tau$ to the sets $\{\f(\U_j)\}_{j=1,\dots,J-1}$ occur in time intervals which do not neighbor each other, i.e. there is a sequence of times $s_1 < t_1 - 1 < s_2 < t_2 - 1 < \dots < s_k < t_k - 1$, and a sequence of indices $(j_1, \dots, j_k)$ in $\{0,\dots, J-1\}$ such that $\tau(t) \in \f(\U_{j_i})$ if $s_i < t < t_i$ ($i = 1,\dots, k$) and $\tau(t) \notin \cup_{j=0}^{J-1}\f(\U_j)$ otherwise.

Using (\ref{eq:skinconn}) and (\ref{eq:tauextr}), we define a first modification $\bar \tau$ of $\tau$ having the same end-points as $\tau$ and which is disjoint from all the `sausages'. We do this by replacing all the pieces $(\tau(t))_{s_i \leqslant t \leqslant t_i}$ by some path in $\partial^* \f(\U_i)$ connecting $\tau(s_1)$ to $\tau(t_1)$. By (\ref{eq:distsauss}) we conclude that $\bar \tau$ is disjoint from $\cup_{j=0}^{J-1} \f(\U_j)$ and using (\ref{eq:diamcondit}), we obtain that $\Range(\bar \tau) \subset B(\Range(\tau),2h)$. Moreover, $\Range(\bar \tau) \subset \Range(\tau) \cup ( \cup_{j = 0}^{J-1} \partial ^* \f(\U_j))$, which by (\ref{eq:distsauss}), is contained in $ \Range(\tau) \cup \partial ^* \widetilde{\U} \cup \{X_{n_u}\}_{j=1,\dots,J-1}$.

In order to find a path $\tau'$ which avoids $\widetilde \U$, we still need to modify $\bar \tau$ in a way that it does not intersect the cut-points $\{X_{n_j} \}$. Recall that the path $\bar \tau$ is disjoint from $\f(\U_j)$, $j=0,\dots,J-1$ and by (\ref{eq:farcuts}), all the neighbors of $X_{n_j}$ which are not in $\f(\U_{j-1}) \cup \f(\U_j)$ are in $Ring_j$.

One can define $\tau'$ by the following: whenever $\bar \tau(t) = X_{n_j}$, the piece $(\bar \tau(t-1), \bar \tau(t), \bar \tau(t+1))$ is replaced by some path in $Ring_j$ connecting $\bar \tau(t-1)$ to $\bar \tau(t+1)$. Since $\diam(Ring_j) = 2$ and $Ring_j \subset \partial^* \widetilde{\U}$, see (\ref{eq:ringskin}) and (\ref{eq:distsauss}), we conclude that $\widetilde{\U}$ is $(3h)$-avoidable.

To finish the proof of the Lemma, we show that
\begin{equation}
\label{eq:fillX}
\widetilde{\U} = \f(X_{(n_0,n_J)}).
\end{equation}

By the Lemma~\ref{eq:charfill}, we conclude that $\f(A) \cup \f(B) \subset \f(A \cup B)$ for any finite sets $A$ and $B$, so that $\widetilde{\U} \subset \f(X_{(n_0,n_J)})$. Since $X_{(n_0,n_J)} \subset \widetilde{\U}$, we only need to prove that $\f(\widetilde{\U}) \subset \widetilde{\U}$. To this end we will show that $\widetilde{\U}^c$ is connected. Given $x, y \in \f(\widetilde{\U})^c$, we connect them by some path $\sigma$. Using the fact that $\widetilde{\U}$ is avoidable, replace this path by some $\tau'$ disjoint from $\widetilde{\U}$, which also joins $x$ to $y$. Thus $\widetilde{\U}^c$ is connected and with Lemma~\ref{eq:charfill} we find that $\widetilde{\U} = \f(\widetilde{\U})$. This shows (\ref{eq:fillX}), finishing the proof of Lemma~\ref{lem:avoid}. $\quad \square$

\vspace*{4mm}

As a consequence of the result above, we prove in the next corollary that, if the diameter of each `sausage' is smaller than $\tfrac{L_0}{4}$, $\f(X_{(n_0,n_J)})$ does not separate components. We write $\chi$ and $C^n$ instead of $\chi_{(0,0)}$ and $C_{(0,0)}^n$ for simplicity.

\begin{corollary}
\label{lem:chifillX}
($d \geqslant 3$) If $n_0 < \dots < n_J$ are cut-points of $(X_n)_{n \in \mathbb{Z}}$ satisfying (\ref{eq:farcuts}) and
\begin{equation}
\label{eq:maxdiam}
\max_{j = 0,\dots,J-1} \{\diam(\U_j) \} < \frac{L_0}{4},
\end{equation}
then $\chi(\f(X_{(n_0,n_J)})=0$ (recall the definition of $\chi$ in (\ref{eq:chi})).
\end{corollary}

\vspace*{4mm}

Note that the cut-times $n_0,\dots, n_J$ are not necessarily all the cut-times in a given interval, see comment under (\ref{eq:Uj}).

\vspace*{4mm}

\textit{Proof.} Take $A_1, A_2 \subset C^2$ such that $\diam(A_i) \geqslant L_0/2$, $i=1,2$. We first show that the conditions on the diameters of the $A_i$'s and $\U_j$'s ensures the existence of a point $x_i$ in $\partial A_i \setminus \f(X_{(n_0,n_j)})$, ($i=1, 2$), c.f. Remark~\ref{rem:Bernoulli} 3).

\vspace*{4mm}

\textit{Case 1: There is some $j \in \{0,\dots,J-1\}$ such that $A_i \cap \partial^*\f(\U_j) \neq \varnothing$.}

In this case, take $y$ in this intersection. If $y = X_{n_j}$ (respectively $X_{n_{j+1}}$) re-choose $y$ as one of its neighbors in $Ring_j$ (resp. $Ring_{j+1}$). In any of these situations, although $y$ can still be an interior point of $A_i$, we know that $y \in \overline{A_i} \setminus \f(X_{(n_0,n_J)})$, see (\ref{eq:fillX}), (\ref{eq:distsauss}) and (\ref{eq:ringskin}). Now, take a path $\tau$ from $y$ to some $y' \in (A_i \cup \f(X_{(n_0,n_J)}))^c$ and use the fact that $\f(X_{(n_0,n_J)})$ is avoidable (see Lemma~\ref{lem:avoid}) to find a $\tau'$ from $y$ to $y'$ which is disjoint from $\f(X_{(n_0,n_J)})$. Finally, we take $x_i$ to be the first point of this path in $\partial A_i$.

\vspace*{4mm}

\textit{Case 2:}
\begin{equation}
\label{eq:AiUfill}
A_i \cap \partial^* \f(\U_j) = \varnothing, \text{ for $j=0, \dots, J-1$}.
\end{equation}
By Claim~\ref{eq:diamfill}, $\diam(\f(\U_j)) = \diam(\U_j) \overset{(\ref{eq:maxdiam})}{<} \diam(A_i)$
implying that $A_i \not \subset \f(\U_j)$ for any $j = 0,\dots, J-1$. This, together with the connectedness of the $A_i$'s and (\ref{eq:AiUfill}) imply that
\begin{equation}
\label{eq:Aidisj}
A_i \cap \overline{\f(\U_j)} = \varnothing, \text{ for $j = 0,\dots, J-1$}.
\end{equation}
In this case, we can choose $x_i$ to be any point in $\partial A_i$ and it will automatically be out of $\f(X_{(n_0,n_J)})$. Otherwise, by (\ref{eq:fillX}), its neighbor in $A_i$ would contradict (\ref{eq:Aidisj}).

\vspace*{4mm}

Now that we have (for $i=1,2$) a point $x_i$ in $\partial A_i\setminus \f(X_{(n_0,n_j)})$, we take any path $\tau$ in $\overline{C^2}$ connecting $x_1$ to $x_2$. Using (\ref{eq:maxdiam}) and the Lemma~\ref{lem:avoid} we obtain a modified path $\tau'$ connecting $\partial A_i$ to $\partial A_2$ which is disjoint from $\f(X_{(n_0,n_J)})$ and contained in $B(C^2,\tfrac{3}{4}L_0+1)$ $\subset C^3$. We proved that $A_1$ and $A_2$ are not separated by $\f(X_{(n_0,n_J)})$ in $C^3$ and since the choice of $A_1$ and $A_2$ was arbitrary, $\chi(\f(X_{(n_0,n_J)}))=0$, concluding the proof of Corollary~\ref{lem:chifillX}. $\quad \square$

\vspace*{4mm}

We now obtain estimates on the the diameter of the `sausages' and on the stopping time $D = T_{C^{L/4}}$ for a typical random walk trajectory, when $d \geqslant 5$.

Recall that our cut-times are defined for doubly infinite trajectories (in fact this has simplified the exposition of Lemma~\ref{lem:avoid}). So, we now artificially introduce a negative time for our random walk trajectory by considering an independent copy of $P_x$. More precisely, let $(X_{-n})_{n \geqslant 0}$ denote the canonical coordinates of the second process on $P_x \otimes P_x$.

\begin{lemma}
\label{lem:estimates}
($d \geqslant 5$) Given $\epsilon > 0$, and integers $G \geqslant 1$ and $L \geqslant 40$, for large enough $L_0 \geqslant 1$ and every $x \in C^5$, with $P_x \otimes P_x$-probability at least $1-\epsilon$ we can find cut times $n_0 < 0 < n_1 < \dots < n_J$ such that:
\begin{equation}
\label{eq:estimates}
\begin{aligned}
i) &\,\, n_0, \dots, n_{J-1} \text{ are not adjacent, see (\ref{eq:farcuts})}, \\
ii)&\,\, {\displaystyle \max_{j=0,\dots,J-1}} \{ \diam(\U_j) \} \leqslant \frac{L_0}{15G}, \\
iii) &\,\, n_J > D.
\end{aligned}
\end{equation}
The number $J$ is deterministic and only depends on $\epsilon$, $G$ and $L$.
\end{lemma}

\vspace*{4mm}

The strategy to prove this lemma can be roughly described as follows. Given integers $F$, $M$ and $K \geqslant 1$, consider the time interval $[-K^2, 2M(F K)^2)$. We split this interval into $2MF^2 + 1$ subintervals of length $K^2$, where we expect to find the cut-times $n_j$. We now make a brief comment on how we are going to pick the constants $M$, $F$ and $K$. Heuristically, we choose:
\begin{equation}
\label{eq:descrestimates}
\begin{aligned}
i) & \,\, \text{$M_0 = M_0(L,\epsilon)$, $K^*(L)$ such that, for every $K \geqslant K^*(L)$, $F \geqslant 1$, with high} \\
& \,\, \text{probability the walk exits $B(0,LF(K+1))$ before the time $M_0(FK)^2$,} \\
ii)& \,\, \text{$F_0 = F_0(G,M_0,\epsilon)$ so that, for every $K \geqslant 1$, with high probability the} \\
& \,\, \text{paths performed by the random walk in each of these time subintervals} \\
& \,\, \text{(of length $K^2$) have diameter at most $\tfrac{F_0K}{90G}$},\\
iii) & \,\, \text{$K_0 = K_0(L,M_0,K^*,F_0,\epsilon) \geqslant K^*$ so that if $K \geqslant K_0$ we can find with high} \\
& \,\, \text{probability at least one cut-time in all these subintervals}.
\end{aligned}
\end{equation}
Finally, given an $L_0 \geqslant F_0 K_0$, we find a $K$ such that $F_0 K \leqslant L_0 \leqslant F_0(K + 1)$ and consider the partition of $[-K^2, 2M_0(F_0 K)^2)$ as above. According to (\ref{eq:descrestimates}) $iii)$, with high probability we can find $2MF^2 + 1$ cut-times, one in each subinterval of our partition. We retain only $J = M_0F_0^2 + 1$ of these cut-times (one in every other subinterval) to ensure that they are not adjacent. Moreover, the choice of constants in (\ref{eq:descrestimates}) $i)$ and $ii)$ will ensure that (\ref{eq:estimates}) $ii)$ and $iii)$ hold with high probability.

\vspace*{4mm}

\textit{Proof of the Lemma~\ref{lem:estimates}.} As we now explain, there is a cut-time in the interval $[0,m)$ with high $P_x \otimes P_x$-probability as $m$ grows, i.e.
\begin{equation}
\label{eq:qmtozero}
q_m \overset{\text{def}}{=} P_x \otimes P_x [\text{there is no cut-time in } [0,m)] \xrightarrow[m \rightarrow \infty]{} 0.
\end{equation}
Recall that our definition of cut-times is slightly different from the definition that for instance appears in \cite{erdos} or \cite{lawler}, p.88. Nevertheless, a slight modification of the argument in \cite{lawler}, p.89 proves that $P_0 \otimes P_0 [ 0 \text{ is a cut-time}] > 0$ (in (\ref{eq:bounddist}) we perform a similar calculation). The statement (\ref{eq:qmtozero}) now follows from the ergodicity of the increments of $X_j$ under $P_0 \otimes P_0$.

Let $T_r$ (for $r > 0$) denote the exit time of the ball $B(0,r)$, recall the definition at the beginning of Section~\ref{sec:review}. We now choose the integers $M_0 \geqslant 1$ and $K^*$, see (\ref{eq:descrestimates}) $i)$. To this end, using the invariance principle, we note that, for large enough $K^*(L) \geqslant 1$,
\begin{equation}
\label{eq:chancetoexit}
\sup_{F \geqslant 1, K \geqslant K^*} P_0 \Big[ \max_{j \leqslant (FK)^2} \lVert X_j \rVert \leqslant 2LF(K+1) \Big] = b < 1.
\end{equation}
Applying the Markov property at the times $(FK)^2,\dots, (M-1)(FK)^2$, we have, for large enough $M_0 = M_0(L,\epsilon)$,
\begin{equation}
\label{eq:exitslow}
P_0[T_{LF(K+1)} > M_0(FK)^2] \leqslant P_0[T_{2LF(K+1)} > (FK)^2]^{M_0} \leqslant b^{M_0} < \epsilon/3,
\end{equation}
for every $K \geqslant K^*$ and $F \geqslant 1$. This completes our choice of $M_0$ and $K^*$ in (\ref{eq:descrestimates}) $i)$.

We now establish estimates on the diameter of the paths performed in each subinterval of length $K^2$, see (\ref{eq:descrestimates}) $ii)$. Let $S_n$ stand for the one-dimensional simple random walk, see below (\ref{eq:departreturn}).
It follows from a variation of Azuma's inequality, see for instance \cite{diarmid}, (41), p.28, that
\begin{equation*}
\begin{split}
P_0 \big[T_{FK/90G} < K^2\big] & \leqslant 2d\, P^1_0\Big[\max_{l \leqslant K^2} S_l > \tfrac{FK}{90G}\Big] \leqslant 4d \, e^{- \tfrac{F^2}{(90G)^2} }, \text{ for every $F$, $K$ and $G \geqslant 1$}.
\end{split}
\end{equation*}
Thus we can choose a large enough $F_0(G,M_0,\epsilon) \geqslant 1$ such that
\begin{equation}
\label{eq:manyexits}
(2M_0F_0^2 + 1) P_0 \big[T_{F_0K/90G} < K^2\big] < \epsilon/3, \text{ for every $K \geqslant 1$}.
\end{equation}

As described in (\ref{eq:descrestimates}) $iii)$, we want to find cut-times in all the $2M_0F_0^2 + 1$ intervals (of length $K^2$) of our partition. For this, using (\ref{eq:qmtozero}), we pick $K_0(L,M_0,K^*,F_0,\epsilon) \geqslant K^*$ so that for every $K \geqslant K_0$,
\begin{equation}
\label{eq:manycuts}
(2M_0F_0^2 + 1)q_{K^2} < \epsilon /3.
\end{equation}

Finally, given an $L_0 \geqslant F_0 K_0$, we choose $K \geqslant K_0$ such that $F_0K \leqslant L_0 < F_0(K+1)$. The bound below is the precise implementation of (\ref{eq:descrestimates}).
\begin{equation}
\label{eq:goodcuttimes}
\begin{split}
P_x \bigg[ & \bigcup_{l = 0}^{2M_0F_0^2} \Big[ \{\text{there is no cut-time in } \big[(l-1)K^2 , lK^2 \big) \} \cup \negthinspace \big\{\diam \big(X_{[(l-1)K^2 , lK^2 )} \big) > \tfrac{F_0K}{45G} \big\} \Big]\\
& \quad \cup \Big\{ T_{LF_0(K + 1)} > M_0(F_0 K)^2 \Big\} \bigg] \overset{(\ref{eq:manycuts})(\ref{eq:manyexits})(\ref{eq:exitslow})}{<} \epsilon.
\end{split}
\end{equation}

On the complement of the event appearing above, we choose $J = M_0F_0^2 + 1$ cut-times in every other time interval (consequently they are not adjacent). For instance, we can choose the first cut-time of the intervals below
\begin{equation*}
n_0 \in \big[-K^2, 0 \big), n_1 \in \big[K^2, 2K^2\big), \dots, n_{M_0F_0^2} \in \big[(2M_0F_0^2-1)K^2, 2M_0(F_0K)^2\big).
\end{equation*}
This ensures (\ref{eq:estimates}) $i)$.

Recall the definition of the sets $U_j$ and the comment below (\ref{eq:Uj}). Since $F_0K \leqslant L_0 < F_0(K+1)$ and $\diam(U_j) \leqslant 3\tfrac{F_0K}{90G}$ on the event in (\ref{eq:goodcuttimes}), we have that (\ref{eq:estimates}) $ii)$ and $iii)$ hold. $\quad \square$

\vspace*{4mm}

\begin{remark}
\label{rem:d3}
\textnormal{There are results concerning ``monolateral'' cut-times of random walks when $d = 3, 4$. In particular, the number of cut-times between zero and $n$ (with a different definition) grows sub-linearly. See for instance, \cite{lawler2}.}
\end{remark}

\vspace*{4mm}

The next result concludes what we called the first step of the proof of Theorem~\ref{th:local}. It establishes that the probability that one random walk trajectory separates macroscopic components inside a box goes to zero as the diameter of the box grows. Moreover, this limit is uniform in the points where we condition the random walk to enter and exit large neighborhoods of the box. Recall that a similar uniformity was important to obtain (\ref{eq:boundsecond}). The proof of the following lemma combines Lemmas~\ref{lem:avoid} and \ref{lem:estimates}, Corollary~\ref{lem:chifillX} and the Harnack inequality (\ref{eq:Harnack}).

Consider the random time
\begin{equation}
\label{eq:timeS}
S \text{ is the last visit of $X$ to $C^4$ before $D$ ($= T_{C^{L/4}}$)}.
\end{equation}
Note that $S$ is not a stopping time.

\begin{lemma}
\label{lem:conditioned}
($d \geqslant 5$) Given $\epsilon >0$, $G \geqslant 1$ and $L \geqslant 40$, for large enough $L_0 \geqslant 1$,
\begin{equation}
\label{eq:conditioned}
\inf_{\substack{x \in \partialint C^4\\ y\in \partial C^{L/4}}} P_{x,y}\Big[ \chi(\f(X_{[0,S]})) = 0, \f(X_{[0,S]}) \text{ is $(\tfrac{L_0}{2G})$-avoidable in $C^3$} \Big] > 1-\epsilon.
\end{equation}
See (\ref{eq:chi}),(\ref{eq:Pxy}) and Definition~\ref{def:avoidable} for the notation.
\end{lemma}

\textit{Proof}. As an intermediate step, we show that with $G$ and $L$ as above, for every $\epsilon' > 0$, there is a large enough $L_0$ such that
\begin{equation}
\label{prop:freewalk}
\inf_{x \in \partialint C^4} P_x
\left[
\begin{array}{c}
\textnormal{for every } 0 \leqslant t < D \textnormal{ such that } X_t \in \partialint C^4, \textnormal{ we have }\\
\chi(\f(X_{[0,t]})) = 0 \textnormal{ and } \f(X_{[0,t]}) \textnormal{ is } \big( \frac{L_0}{2G} \big)\textnormal{-avoidable in } C^3
\end{array}
\right] > 1-\epsilon'.
\end{equation}

To this end, we first take $L_0$ as in Lemma~\ref{lem:estimates}. We know that on an event with probability at least $1-\epsilon'$ we can choose the cut-times $n_0 < 0 < \dots < n_J$ satisfying (\ref{eq:estimates}), recall the observation below (\ref{eq:Uj}). On this event, for any $0 \leqslant t < D$ such that $X_t \in \partialint C^4$, we take $\bar j$ such that $n_{\bar j} \leqslant t < n_{\bar j +1}$.

Using (\ref{eq:chimonotone}) and Corollary~\ref{lem:chifillX}, we conclude that
\begin{equation}
\chi \big(\f(X_{[0,t]}) \big) \leqslant \chi \big(\f(X_{(n_0,n_{\bar j +1})}) \big) = 0.
\end{equation}

Since (\ref{eq:estimates}) $ii)$ holds, we can take $h = \tfrac{L_0}{6G} > \tfrac{L_0}{15G}$ in Lemma~\ref{lem:avoid} to obtain that $\f(X_{(n_0,n_{\bar j +1})})$ is $(\tfrac{L_0}{2G})$-avoidable in $\mathbb{Z}^d$ and by (\ref{eq:avoidmonotone}),
\begin{equation}
\label{eq:X0j1avoidable}
\f(X_{(n_0,n_{\bar j +1})}) \text{ is $(\tfrac{L_0}{2G})$-avoidable in $C^3$}.
\end{equation}
Since
\begin{equation}
\label{eq:X0tsurrounded}
\f(X_{(n_1,n_{\bar j})}) \subset \f(X_{[0,t]}) \subset \f(X_{(n_0,n_{\bar j +1})}),
\end{equation}
according to (\ref{eq:avoidremark}) and (\ref{eq:X0j1avoidable}), all we need in order to show that $\f(X_{[0,t]})$ is also $(\tfrac{L_0}{2G})$-avoidable in $C^3$ is that
\begin{equation}
\label{eq:samefill}
\f(X_{(n_1,n_{\bar j})}) \cap B(C^3,\tfrac{L_0}{2G}) = \f(X_{(n_0,n_{\bar j + 1})}) \cap B(C^3,\tfrac{L_0}{2G}).
\end{equation}
Using (\ref{eq:fillX}), we obtain that $\f(X_{(n_0,n_{\bar j +1})}) \setminus \f(X_{(n_1,n_{\bar j})})$ is contained in $\overline{\f(\U_0) \cup \f(\U_{n_{\bar j}})}$. But since $X_0, X_t \in \partialint C^4$ and $n_0 \leqslant 0 < n_1, n_{\bar j} \leqslant t < n_{\bar j +1}$, we know by (\ref{eq:estimates}) $ii)$ that $\overline{\f(\U_0) \cup \f(\U_{n_{\bar j}})}$ is disjoint from $B(C^3,\tfrac{L_0}{2G})$. This establishes (\ref{eq:samefill}) and consequently (\ref{prop:freewalk}).

We now introduce the stopping time
\begin{equation}
\begin{split}
S' = \inf \big\{s \geqslant 0; & X_s \in \partialint C^4, \chi(\f(X_{[0,s]})) = 1 \text{ or}\\
& \f(X_{[0,s]}) \text{ is not $(\tfrac{L_0}{2G})$-avoidable in $C^3$} \big\},
\end{split}
\end{equation}
and note that the event appearing in (\ref{eq:conditioned}) contains $\{S' \geqslant D\}$.

Define, for $y \in \partial C^{L/4}$ and $z \in C^{L/4}$ the function $h^y(z) = P_z[X_D = y]$ which is harmonic in $C^{L/4}$. Given $x \in \partialint C^4$ and $y \in \partial C^{L/4}$, we use the strong Markov property at time $S'$ to obtain
\begin{equation}
P_{x,y}[S'<D] = \frac{P_x[S'<D,X_D = y]}{h^y(x)} \leqslant \sup_{x,z \in C^4} \frac{h^y(z)}{h^y(x)} P_x[S'<D].
\end{equation}

By the Harnack inequality, see \cite{lawler} Theorem~1.7.6 p.46, we have
\begin{equation}
\label{eq:Harnack}
\sup_{L_0 \geqslant 1} \quad \sup_{x,z \in C^5} \frac{h^y(z)}{h^y(x)} = c < \infty,
\end{equation}
and the Lemma~\ref{lem:conditioned} follows from (\ref{prop:freewalk}) by choosing $\epsilon' = \epsilon/c$. $\quad \square$

\vspace*{4mm}

In Theorem~\ref{th:local}, one considers $G$ independent paths instead of just one as in the above lemma, see also Remark~\ref{rem:G100}. The following lemma is the key step to obtain this extension.

As explained at the beginning of this section, Lemma~\ref{lem:farapart} shows that for any family of $G$ connected sets, which do not separate components, are avoidable and mutually far apart, the union of these sets also does not separate components. More precisely,


\begin{lemma}
\label{lem:farapart}
Let $U_1,\dots, U_G \subset \mathbb{Z}^d$ be connected sets. Setting $\mathcal{O}_1 = \f(U_1), \dots, \mathcal{O}_G = \f(U_G) \subset \mathbb{Z}^d$, if the following holds:
\begin{equation}
\label{eq:chiOi}
\begin{aligned}
i) \,\,& \chi(\mathcal{O}_i) = 0, \text{ for all } i=1,\dots, G,\\
ii) \,\,& d_\infty(U_i,U_j) > 1 \text{ for all } 1 \leqslant i < j \leqslant G,\\
iii)\,\,& \mathcal{O}_i \text{ is $(\tfrac{L_0}{2G})$-avoidable in $C^3$ for every } i=1,\dots,G,
\end{aligned}
\end{equation}
then, $\chi(\cup_{i=1}^G \mathcal{O}_i) = 0.$
\end{lemma}

\vspace*{4mm}

\textit{Proof.} As we now show, we can assume without loss of generality that
\begin{equation}
\label{eq:redfarapart}
\text{$\mathcal{O}_i \setminus \mathcal{O}_j \neq \varnothing$ for all distinct $i,j \in \{1,\dots, G\}$}.
\end{equation}

Indeed, if $\mathcal{O}_i \subset \mathcal{O}_j$ for some distinct pair $1 \leqslant i,j \leqslant G$, we eliminate this $\mathcal{O}_i$. So we can assume (\ref{eq:redfarapart}). Hence, with Lemma~\ref{prop:fill} and (\ref{eq:chiOi}) $ii)$ we conclude that
\begin{equation}
\label{eq:dOiOj}
d_\infty(\mathcal{O}_i,\mathcal{O}_j) > 1 \text{ for all } 1 \leqslant i < j \leqslant G.
\end{equation}

To obtain $\chi(\cup_{i=1}^G \mathcal{O}_i) = 0$ (see (\ref{eq:chi})), we will prove that
\begin{equation}
\label{eq:dA1A2anddiam}
\begin{array}{c}
\text{for every pair of sets $A_1, A_2 \subset C^2$ such that $A_1$ and $A_2$,}\\
\text{are connected, $d(A_1,A_2) > 1$ and $\diam(A_1), \diam(A_2) \geqslant L_0/2$ }\\
\text{there is a path in $C^3 \setminus \cup_{i=1}^G \mathcal{O}_i$ connecting $\partial A_1$ to $\partial A_2$}
\end{array}
\end{equation}

As a further reduction, we are going to prove that (\ref{eq:dA1A2anddiam}) follows if one shows (\ref{eq:dA1A2anddiam}) when
\begin{equation}
\label{eq:redfarapart2}
\begin{array}{c}
\text{$A_1$ and $A_2$ in (\ref{eq:dA1A2anddiam}) satisfy the additional condition $d(\f(A_1), \f(A_2)) > 1$}.
\end{array}
\end{equation}

To prove the above reduction, it suffices to show the following fact. Given any pair $A_1, A_2 \subset \mathbb{Z}^d$ satisfying (\ref{eq:dA1A2anddiam}), but such that $d(\f(A_1), \f(A_2)) \leqslant 1$,
\begin{equation}
\label{eq:pathtouches}
\begin{array}{c}
\text{we can exhibit sets $A_1'$ and $A_2'$ as in (\ref{eq:dA1A2anddiam}) and fulfilling the additional} \\
\text{condition in (\ref{eq:redfarapart2}) such that: if there is no path in $C^3 \setminus \cup_{i=1}^G \mathcal{O}_i$ connecting} \\
\text{$\partial A_1$ to $\partial A_2$, then there is also no such path between $\partial A_1'$ and $\partial A_2'$.}
\end{array}
\end{equation}
In other words, if $\cup_{i=1}^G \mathcal{O}_i$ separates $A_1$ from $A_2$ in $C^3$, it also separates $A_1'$ from $A_2'$ in $C^3$.

\begin{figure}[ht]
\psfrag{a1}{$A_1$}
\psfrag{a2}{$A_2$}
\psfrag{F}{$A_1'$}
\begin{center}
\includegraphics[angle=0, width=0.25\textwidth]{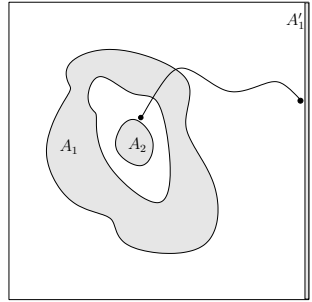}\\
\caption{Every path joining $\partial A_2' = \partial A_2$ to $\partial A_1'$ meets $\partial A_1$.}\label{fig:A1A2F}
\end{center}
\end{figure}

Let us now explain how this fact is proved. Using (\ref{eq:dA1A2anddiam}) and Lemma~\ref{prop:fill}, since $d(\f(A_1), \f(A_2)) \leqslant 1$, we know that either $A_2 \subset \f(A_1)$ or $A_1 \subset \f(A_2)$. We suppose without loss of generality that we are in the former case. We choose the sets $A_2' = A_2$ and $A_1'$ to be some face of $\partialint C^2$, see Figure~\ref{fig:A1A2F}. Since $A_2 = A_2'$, to establish (\ref{eq:pathtouches}), is enough to show that any path connecting $\partial A_2'$ to $\partial A_1'$ must also intersect $\partial A_1$. This is done in the next paragraph.

Suppose by contradiction that there is a path from $\partial A_2'$ to some $y \in \partial A_1'$, which does not meet $\partial A_1$ (and which by (\ref{eq:dA1A2anddiam}) does not meet $A_1$), we could continue this path to the neighbor of $y$ in $A_1'$, then to $(C^2)^c$ and finally to infinity, without touching $A_1$. This would contradict, using Lemma~\ref{eq:charfill}, the fact that $A_2 \subset \f(A_1)$. Using a similar, although simpler reasoning, we also obtain that $A_1' \cap \overline{A_2'} = \varnothing$, implying that $d(\f(A_1'),A_2') = d(A_1',A_2') > 1$. Since $\diam(A_1') \geqslant L_0/2$, we proved (\ref{eq:pathtouches}) and consequently (\ref{eq:redfarapart2}).

Given a pair $A_1, A_2 \subset \mathbb{Z}^d$, satisfying the conditions in (\ref{eq:dA1A2anddiam}) and (\ref{eq:redfarapart2}), we will exhibit a path in $C^3 \setminus \cup_{i=1}^G \mathcal{O}_i$ connecting $\partial A_1$ to $\partial A_2$. First we find a point $x_1$ in $\partial A_1 \setminus \cup_{i=1}^G \mathcal{O}_i$. For this, note that by (\ref{eq:fillstconn}),
\begin{equation}
\label{eq:pfillA1conn}
\partial \f(A_1) \text{ is $*$-connected}
\end{equation}
and by (\ref{eq:redfarapart2}), we can assume $d(\f(A_1), \f(A_2)) > 1$, so that
\begin{equation}
\label{eq:pfillA1separates}
\partial \f(A_1) \text{ separates $A_1$ from $A_2$ in $C^3$}.
\end{equation} 

In order to find a point $x_1$ in $\partial A_1 \setminus \cup_{i=1}^G \mathcal{O}_i$, we first take any $x' \in \partial \f(A_1)$, which by (\ref{eq:boundfill}) is also in $\partial A_1$. If $x' \not \in \cup_{i=1}^G \mathcal{O}_i$, we are done, otherwise, let $i_o$ be such that $x' \in \mathcal{O}_{i_o}$. By (\ref{eq:chiOi}) $i)$ and (\ref{eq:pfillA1separates}), we are able to find some $x'' \in \partial \f(A_1) \setminus \mathcal{O}_{i_o}$. Using the $*$-connectedness of $\partial \f(A_1)$ (see (\ref{eq:pfillA1conn})) we join $x'$ to $x''$ by a $*$-path $\sigma$ in $\partial \f(A_1)$, and take $x_1$ to be the first point of $\sigma$ out of $\mathcal{O}_{i_o}$. We conclude from (\ref{eq:dOiOj}) that $x_1 \in \partial A_1 \setminus \cup_{i=1}^G \mathcal{O}_i$.

In the same way, we find some $x_2 \in \partial A_2 \setminus \cup_{i=1}^G \mathcal{O}_i$ and join $x_1$ to $x_2$ by any path $\tau \subset C^2$. Roughly speaking, to conclude the proof we will modify $G$ times the path $\tau$ (using (\ref{eq:chiOi}) $iii)$) in order to avoid each set $\{\mathcal{O}_i\}_{i=1,\dots,G}$.

Since $\mathcal{O}_1$ is $(\tfrac{L_0}{2G})$-avoidable in $C^3$ (see (\ref{eq:chiOi}) $iii)$), we can find a modification $\tau_1$ of $\tau$, joining $x_1$ to $x_2$, which is disjoint from $\mathcal{O}_1$ and such that $\Range(\tau_1) \subset B(C^2, \tfrac{L_0}{2G})$.

We proceed by induction. Suppose that for some $1 \leqslant j < G$ we found some $\tau_j$ joining $x_1$ to $x_2$ such that
\begin{equation}
\begin{aligned}
i) & \,\, \text{$\Range(\tau_j) \cap \mathcal{O}_i = \varnothing$ for $i\leqslant j$ and}\\
ii) & \,\, \text{$\Range(\tau_j) \subset B(C^2, \tfrac{jL_0}{2G})$}.
\end{aligned}
\end{equation}
We use (\ref{eq:chiOi}) $iii)$, see also Definition~\ref{def:avoidable}, to find a path $\tau_{j+1}$ joining $x_1$ to $x_2$ such that
\begin{equation}
\begin{aligned}
i) & \,\, \text{$\Range(\tau_{j+1}) \subset (\Range(\tau_j) \setminus \mathcal{O}_{j+1}) \cup \partial^* \mathcal{O}_{j+1}$, implying by (\ref{eq:dOiOj}) that:}\\
& \,\, \text{$\Range(\tau_{j+1}) \cap \mathcal{O}_i = \varnothing$, for $i \leqslant j+1$ and}\\
ii) & \,\, \text{$\Range(\tau_{j+1}) \subset B(C^2,\tfrac{(j+1)L_0}{2G})$}.
\end{aligned}
\end{equation}

The existence of $\tau_G$ as above implies $\chi(\cup_{i=1}^G \mathcal{O}_i) = 0$ and consequently, Lemma~\ref{lem:farapart}. $\quad \square$

\vspace*{4mm}

The next Lemma is the final ingredient to prove the main result of this section. It will ensure that with high probability a set of $G$ independent random walks satisfy the hypothesis (\ref{eq:chiOi}) $ii)$ of the Lemma~\ref{lem:farapart}, or in other words: they are mutually far apart.

The proof of this lemma is an adaptation of known arguments concerning intersection of random walk trajectories for $d \geqslant 5$, see for instance, \cite{lawler} p.89.

We denote by $D^i$, $H^i_K$ and $S^i$ the times $D$, $H_K$ and $S$ pertaining to the walks $(X^i)_{i=1,\dots, G}$, recall (\ref{eq:timeS}).

\begin{lemma}
\label{lem:Gfarapart}
($d \geqslant 5$) Given $\epsilon > 0$ and $L \geqslant 40$, for large enough $L_0 \geqslant 1$,
\begin{equation}
\begin{split}
\sup_{\substack{x_1 \in \partialint C^4,x_2 \in \partialint C^5 \\ y_1,y_2 \in \partial C^{L/4}}}
P_{x_1,y_1}\otimes P_{x_2,y_2} \Big[H^2_{C^4} < D^2, d_\infty(X_{[0 , S^1]}^1, X_{[H_{C^4}^2, S^2]}^2) \leqslant 1 \Big]< \epsilon
\end{split}
\end{equation} 
\end{lemma}

\vspace*{4mm}

\textit{Proof.} As in the proof of Lemma~\ref{lem:conditioned}, we define $h^y(z) = P_z[X_D = y]$, for $z\in C^{L/4}$, $y\in \partial C^{L/4}$, recall that for $x, y \in \mathbb{Z}^d$ we write $x \overset{*}{\leftrightarrow} y$ if they are $*$-neighbors. For $x_1 \in \partialint C^4$, $x_2 \in \partialint C^5$ and $y_1,y_2 \in \partial C^{L/4}$, we have
\begin{equation}
\begin{array}{l}
\displaystyle
\vspace{2mm}
P_{x_1,y_1}\otimes P_{x_2,y_2} \Big[H^2_{C^4} < D^2,d_\infty(X_{[0 , S^1]}^1, X_{[H_{C^4}^2, S^2]}^2) \leqslant 1 \Big]\\
\vspace{2mm}
\displaystyle \leqslant P_{x_1,y_1}\otimes P_{x_2,y_2} \bigg[
\begin{array}{c}
\text{there are } s_1 \leqslant t_1 < D^1, s_2 \leqslant t_2< D^2;\\
X_{t_1}^1, X_{t_2}^2 \in C^4 \text{ and } X_{s_1}^1 \overset{*}{\leftrightarrow} X_{s_2}^2
\end{array}
\bigg]\\
\displaystyle \leqslant \sum_{s_1 \geqslant 0} \sum_{s_2 \geqslant 0} \,\,\, \sum_{\substack{z_1,z_2 \in C^{L/4}; z_1 \overset{*}{\leftrightarrow} z_2}} P_{x_1,y_1}[s_1 < D,X_{s_1} = z_1] P_{z_1,y_1}[H_{C^4} < D]\\
\displaystyle \phantom{= \sum_{s_1 \geqslant 0} \sum_{s_2 \geqslant 0} \sum_{z_1 \overset{*}{\leftrightarrow} z_2}} \quad P_{x_2,y_2}[s_2 < D,X_{s_2} = z_2] P_{z_2,y_2}[H_{C^4} < D]\\
\displaystyle \leqslant \sum_{s_1 \geqslant 0} \sum_{s_2 \geqslant 0} \,\,\, \sum_{z_1,z_2 \in C^{L/4}; z_1 \overset{*}{\leftrightarrow} z_2} P_{x_1}[s_1 < D,X_{s_1} = z_1] \frac{h^y(z_1)}{h^y(x_1)} \sup_{w_1 \in C^4} \frac{h^y(w_1)}{h^y(z_1)}\\
\displaystyle \phantom{= \sum_{s_1 \geqslant 0} \sum_{s_2 \geqslant 0} \sum_{z_1 \overset{*}{\leftrightarrow} z_2}} \quad P_{x_2}[s_2 < D,X_{s_2} = z_2] \frac{h^y(z_2)}{h^y(x_2)} \sup_{w_2 \in C^4} \frac{h^y(w_2)}{h^y(z_2)},
\end{array}
\end{equation}
using reversibility and Harnack's inequality (\ref{eq:Harnack}) for $\displaystyle \sup_{w_i \in C^4} \frac{h^y(w_i)}{h^y(x_i)}$, $i = 1,2$, (note that we cannot use it with the $z_i$'s since they could be out of $C^5$),
\begin{equation}
\begin{array}{l}
\displaystyle \leqslant c^2 \sum_{s_1 \geqslant 0} \sum_{s_2 \geqslant 0} \sum_{z_1,z_2 \in \mathbb{Z}^d;  z_1 \overset{*}{\leftrightarrow} z_2} P_{x_1}[s_1 < D,X_{s_1} = z_1] P_{z_2}[s_2 < D,X_{s_2} = x_2],
\end{array}
\end{equation}
linking $z_1$ to $z_2$ in at most $d$ steps, we bound the expression above by
\begin{equation}
\begin{array}{l}
(2d)^dc^2 \sum_{r \geqslant 0} (r+1) P_{x_1}[X_r = x_2].
\end{array}
\end{equation}
Finally we use the heat kernel estimate $P_x[X_n = y] \leqslant cn^{d/2}\exp(\tfrac{|x-y|^2}{cn})$, see \cite{grigoryan} (2.4), to obtain
\begin{equation}
\label{eq:bounddist}
P_{x_1,y_1}\otimes P_{x_2,y_2} \Big[H^2_{C^4} < D^2,d_\infty(X_{[0 , S^1]}^1, X_{[H_{C^4}^2, S^2]}^2) \leqslant 1 \Big] \leqslant c \sum_{r \geqslant 0} \frac{r+1}{r^{d/2}} e^{-c'\tfrac{L_0^2}{r}}.
\end{equation}
This last quantity, independently on the choice of $x_1$, $x_2$, $y_1$, $y_2$, goes to zero as $L_0$ goes to infinity (recall that $d \geqslant 5$), and Lemma~\ref{lem:Gfarapart} follows. $\quad \square$

\vspace*{4mm}

The next theorem is the main result of this section, recall that it was already stated and used in Section~\ref{sec:proof}. Its proof combines Lemmas~\ref{lem:conditioned}, \ref{lem:farapart} and \ref{lem:Gfarapart}. Recall the definition of the time $S$ in (\ref{eq:timeS}) and the sets $\Delta_1$, $\Delta_2$ below (\ref{eq:Pxy}).

\begin{theorem}
\label{th:local}
($d \geqslant 5$) Given $\epsilon > 0$, $G \geqslant 1$ and $L \geqslant 40$, for large enough $L_0 \geqslant 1$,
\begin{equation}
\label{eq:localbound}
\sup_{\mbox{\fontsize{8}{8} \selectfont $\substack{\vec{x} \in \Delta_1 \\ \vec{y} \in \Delta_2}$}} P_{\vec{x},\vec{y}} \left [\chi \bigg( \bigcup_{i=1}^G X^i_{[0,D]} \bigg) = 1 \right] < \epsilon.
\end{equation}
\end{theorem}

\textit{Proof.} Fix $\vec{x} \in \Delta_1$ and $\vec y \in \Delta_2$ and recall that for any $\U \subset \mathbb{Z}^d$, $\chi(\U)$ depends only on $\U \cap C^3$.

Roughly speaking, in Lemma~\ref{lem:farapart} we have seen that if a union of $G$ sets separates components in $C^3$ then: these sets are not mutually far apart, or one of them either separates components in $C^3$ or is not $(\tfrac{L_0}{2G})$-avoidable in $C^3$. More precisely, Lemma~\ref{lem:Gfarapart} yields the following:
\begin{equation*}
\begin{array}{l}
\vspace*{2mm}
{\displaystyle P_{\vec{x},\vec{y}}\bigg[ \chi \Big(\bigcup_{i=1}^{G} X_{[0,D]}^i \Big) = 1 \bigg]}
\leqslant P_{\vec{x},\vec{y}}
\left[
\begin{array}{c}
\text{there are } 1 \leqslant i < j \leqslant G \text{ such that $H_{C^4}^i < D^i$,}\\
H_{C^4}^j < D^j \text{ and } d_\infty(X_{[H_{C^4}^i,S^i]}^i,X_{[H_{C^4}^j,S^j]}^j) \leqslant 1 
\end{array}
\right]
\end{array}
\end{equation*}
\begin{equation*}
\begin{array}{l}
\vspace*{2mm}
\quad + P_{\vec{x},\vec{y}}
\left[
\begin{array}{c}
\text{there is some } 1 \leqslant i \leqslant G \text{ such that } H_{C^4}^i < D^i \text{ and either }\\
\chi(\f(X_{[H_{C^4}^i,S^i]}^i)) = 1 \text{ or } \f(X_{[H_{C^4}^i,S^i]}^i) \text{ is not $(\tfrac{L_0}{2G})$-avoidable in $C^3$}
\end{array}
\right],
\end{array}
\end{equation*}
and using in both terms the strong Markov property for $X^i$ at time $H^i_{C^4}$,
\begin{equation*}
\begin{array}{l}
\leqslant  {\displaystyle \frac{G(G-1)}{2} \sup_{\substack{x_1 \in \partialint C^4, x_2 \in \partialint C^5 \\ y_1,y_2 \in \partial C^{L/4}}}}
P_{x_1,y_1}\otimes P_{x_2,y_2} \Big[H^2_{C^4}< D^2,d_\infty(X_{[0 , S^1]}^1, X_{[H_{C^4}^2, S^2]}^2) \leqslant 1 \Big]\\
\end{array}
\end{equation*}
\begin{equation*}
\begin{array}{l}
\vspace*{2mm}
+ {\displaystyle \, G \sup_{\substack{x \in \partialint C^4\\ y\in \partial C^{L/4}}}} P_{x,y} \Big[ \chi(\f(X_{[0,S]})) = 1 \text{ or } \f(X_{[0,S]}) \text{ is not $(\tfrac{L_0}{2G})$-avoidable in $C^3$} \Big],
\end{array}
\end{equation*}
which, by the Lemmas~\ref{lem:Gfarapart} and \ref{lem:conditioned}, can be made arbitrarily small once we choose $L_0$ large enough. Finishing the proof of the Theorem~\ref{th:local}. $\quad \square$

\vspace*{4mm}

This concludes the proof of (\ref{th:reftolocal}) and provides the last missing piece of the proof of our main Theorem~\ref{th:main}.

\begin{remark}
\textnormal{The present work leaves several questions untouched. For instance:}

\textnormal{- Can one improve Theorem~\ref{th:taildiam} in such a way that the exponents of $N$ in the lower and upper bound match? If the answer is affirmative, what is this exponent? The same questions can be asked for $V$ in Theorem~\ref{th:tailvol}.}

\textnormal{- Which results of this article can be extended to any $u < u_*$ or for $d = 3, 4$? See the Remark~\ref{rem:Bernoulli} 4), the beginning of Section~\ref{sec:local} and Remark~\ref{rem:d3}.}

\textnormal{- How does the size of $\mathcal{C}^u_0$ behave in the sub-critical phase $u > u_*$?}
$\quad \square$.
\end{remark}

\appendix

\section*{Appendix: Properties of $\f(A)$}
\label{app:fill}

In this appendix we prove Lemma~\ref{prop:fill} which we used often in Section~\ref{sec:local}.

\vspace*{4mm}
\textit {Proof of Lemma~\ref{prop:fill}.} First we prove that $\partialint \f(A) \subset A$, see (\ref{eq:boundfill}). If $z \in \partialint \f(A) \setminus A$, one can join $z$ to $\f(A)^c$ and then to infinity without meeting $A$, a contradiction to the assumption that $z \in \f(A)$, in view of Lemma~\ref{eq:charfill}. This proves (\ref{eq:boundfill}).

To show that $\diam(\f(A)) = \diam(A)$, recall (\ref{eq:diamfill}), it is enough to see that
\begin{equation}
\diam(A) \leqslant \diam(\f(A)) \overset{Lemma~\ref{eq:charfill}}{\leqslant} \diam (\text{sbox}(A)) = \diam(A).
\end{equation}

Finally we prove (\ref{eq:threecases}), the last part of Lemma~\ref{prop:fill}. Consider finite and connected sets $A, B \subset \mathbb{Z}^d$. We claim that, if $d(A,B) > 1$, then $B \cap \partial (\f(A)) = A \cap \partial (\f(B)) = \varnothing$. Indeed, by Claim~\ref{eq:boundfill}, every point of $\partial \f(A)$ is a neighbor of $A$, so that they cannot belong to $B$. The same argument applies replacing $A$ with $B$.

Since $B$ is connected and does not intersect $\partial \f(A)$, we have either $B \subset \f(A)$ or $B \cap \overline{\f(A)} = \varnothing$. Analogously, we obtain $A \subset \f(B)$ or $A \cap \overline{\f(B)} = \varnothing$.

It is not possible that $B \subset \f(A)$ and $A \subset \f(B)$. This would imply $\f(A) = \f(B)$ and by Claim~\ref{eq:boundfill} that $d(A,B) = 0$. So we have three remaining possibilities which correspond to the cases enumerated in Lemma~\ref{prop:fill}.

Our claim will now follow once we show that $B \cap \overline{\f(A)} = \varnothing$ and $A \cap \overline{\f(B)} = \varnothing$ imply that $d(\f(A),\f(B)) > 1$ (respectively $d_\infty(\f(A),\f(B)) > 1$).

First we show that $\f(A) \cap \f(B) = \varnothing$. Since $\f(A) \cap \f(B)$ is finite, it is enough to show that the boundary of this set is empty. Suppose $x \in \f(A) \cap \f(B)$ and let $y$ be a neighbor of $x$. If $y \notin \f(A)$, then by Claim~\ref{eq:boundfill} $x \in A$ in contradiction with the fact that $A \cap \overline{\f(B)} = \varnothing$. Using a symmetric argument for $\f(B)$, we conclude that $y \in \f(A) \cap \f(B)$, showing that the boundary of this set is empty. This implies that $\f(A) \cap \f(B) = \varnothing$.

We then prove that $d(\f(A),\f(B)) > 1$ by excluding the possibility that there are neighbors $x \in \f(A)$ and $y \in \f(B)$. Suppose by contradiction the existence of such $x$ and $y$. Since $\f(A) \cap \f(B) = \varnothing$, $y \notin \f(A)$, implying that $x \in \partialint \f(A) \subset A$, by Claim~\ref{eq:boundfill}, a contradiction with $d(\f(B),A) > 1$.

If in addition $d_\infty(A,B) > 1$, we now show that $d_\infty(\f(A),\f(B)) > 1$. Indeed, suppose by contradiction that $x \in \f(A)$ and $y \in \f(B)$ are $*$-neighbors. Since $\f(A) \cap \f(B) = \varnothing$, we conclude that $x \not \in \f(B)$ and $y \not \in \f(A)$. Take a path $\tau$ joining $x$ to $y$ such that $\diam(\Range(\tau)) = 1$. We denote by $x'$ the last visit of $\tau$ to $\f(A)$ and by $y'$ the first visit of $\tau$ to $\f(B)$. Since $x' \in \partialint \f(A) \subset A$ and $y' \in \partialint \f(B) \subset B$, we obtain a contradiction with the hypothesis $d_\infty(A,B) > 1$. $\quad \square$

\end{document}